\documentclass[11pt]{article}
\usepackage{amsmath,amsfonts,graphics,color,amssymb,amstext}
\newtheorem{thm}{Theorem}[section]
\newtheorem{rem}[thm]{Remark}
\newtheorem{lemma}[thm]{Lemma}

\newtheorem{cor}[thm]{Corollary}
\newtheorem{defn}[thm]{Definition}
\newenvironment{pf}{\noindent{\sc Proof. }}{{\hfill$\Box$}\par\vskip2\parsep}

\newcommand{\Let}{Let\\[-2em]}
\newcommand{\llet}{let\\[-2em]}

\newcommand{\ey}{\mbox{\bf E}}
\newcommand{\eyb}{\mbox{\bf \Large E}}

\newcommand{\N}{{\mathbb N}}
\newcommand{\Z}{{\mathbb Z}}

\newcommand{\odd}{\mathrm{od}}
\newcommand{\even}{\mathrm{ev}}
\newcommand{\zd}{{{\mathbb Z}^d}}
\newcommand{\ed}{{{\mathbb Z}^d_\even }}
\newcommand{\od}{{{\mathbb Z}^d_\odd }}

\newcommand{\edn}{{\mathcal E}^d_{N}}

\newcommand{\odl}{{\Lambda_{\odd} }}

\newcommand{\tnod}{{\mathbb T}^d_{N,\odd}}
\newcommand{\EEN}{{\mathbb T}^d_{N,\even}}
\newcommand{\tnev}{{\mathbb T}^d_{N,\even}}

\newcommand{\F}{{\mathcal F}}
\newcommand{\A}{{\mathcal A}}
\newcommand{\C}{{\mathcal C}}

\newcommand{\bzd}{\ensuremath{{(\zd)^*}}}

\newcommand{\bed}{\ensuremath{{(\ed)^*}}}
\newcommand{\bzdl}{\ensuremath{{\Lambda^*}}}
\newcommand{\bzdlb}{\ensuremath{{\overline{\Lambda^*}}}}

\newcommand{\bedl}{\ensuremath{{(\Lambda_\even)^*}}}

\newcounter{mycount}

\newlength{\slimlen}%
\newcommand{\ormd}{\mathrm{d}}
\newcommand{\rmd}{\,\mathrm{d}}

\newcommand{\ZZ}{\ensuremath{\mathbb{Z}}}
\newcommand{\TT}{\ensuremath{\mathbb{T}}}

\newcommand{\RR}{\ensuremath{\mathbb{R}}}
\newcommand{\btd}{\nabla}
\newcommand{\var}{{\mbox{\rm\normalfont Var}}}

\newcommand{\ext}{{\mbox{\rm\normalfont ext\;}}}

\newcommand{\cov}{{\mbox{\rm\normalfont cov}}}

\usepackage{graphics}

\title{Decay of covariances, uniqueness of ergodic component and scaling limit for a class of $\nabla\phi$ systems with non-convex potential}

\author{
Codina Cotar
\thanks{Supported by the DFG-Forschergruppe 718 \lq Analysis and stochastics in
complex physical systems\rq}
\thanks{Corresponding Author}
\thanks{%
TU M\"unchen - Zentrum Mathematik,
Lehrstuhl f\"ur Mathematische Statistik  ,
Boltzmannstr. 3, 85747 Garching, Germany.
E-mail: {\tt cotar@ma.tum.de}}
  and Jean-Dominique Deuschel\footnotemark[1]\;\thanks{%
TU Berlin - Fakult\"at II
Institut f\"ur Mathematik
Strasse des 17. Juni 136
D-10623 Berlin, Germany.
E-mail: {\tt deuschel@math.tu-berlin.de}}}

\begin{document}
\maketitle
\begin{abstract}
We consider a gradient interface model on the lattice with interaction potential which is a
non-convex perturbation of a convex potential. Using a technique which decouples the neighboring vertices into even and odd vertices, we show for a class of non-convex potentials: the uniqueness of ergodic component for $\nabla\phi$- Gibbs measures, the decay of covariances, the scaling limit and the strict convexity of the surface tension.
\end{abstract}


         

\thispagestyle{empty}

\par
---------------------------------
\section{Introduction}

\subsection{The setup}

Phase separation in $\RR^{d+1}$ can be described by effective interface models, where interfaces are sharp boundaries which separate the different regions of space occupied by different phases. In this class of models, the interface is modeled as the graph of a random
function from $\Z^d$ to $\Z$ or $\RR$ (discrete or continuous effective interface models). For more on interface models, see the reviews by Funaki \cite{FSL} or Velenik \cite{vel}. In this setting we ignore overhangs and for $x\in\zd$, we denote by $\phi(x)\in\RR$
the height of the interface above or below the site $x$.
Let $\Lambda$ be a finite set in $\zd$ with boundary 
\begin{eqnarray}
\partial\Lambda:=\{x\notin\Lambda,~||x-y||=1~\mbox{for some}~y\in\Lambda\},~\mbox{where}~\|x-y\|_=\sum^d_{i=1}|x_i-y_i|~\mbox{for}~x,y\in\Z^d
\end{eqnarray}
and with given boundary condition $\psi$ such that $\phi(x)=\psi(x)$ for $x\in\partial\Lambda$; a special case of boundary conditions are the {\it tilted} boundary conditions, with $\psi(x)=x\cdot u$ for all $x\in\partial\Lambda$, and where $u\in\RR^d$ is fixed. Let $\bar\Lambda:=\Lambda\cup\partial\Lambda$ and let $\rmd\phi_{\Lambda}=\prod_{x\in\Lambda}\rmd\phi(x)$ be the Lebesgue measure over $\RR^{\Lambda}$.
For a finite region $\Lambda\subset\zd$, {\it the finite volume Gibbs measure $\nu_{\Lambda,\psi}$ on $\RR^{\zd}$ with boundary condition $\psi$ }for the field of \textit{height variables} $(\phi(x))_{x\in\zd}$ over $\Lambda$ is defined by
\begin{equation}
\label{tag0'}
\nu_{\Lambda,\psi}(\ormd\phi)=\frac{1}{Z_{\Lambda,\psi}}\exp\left\{-\beta H_{\Lambda,\psi}(\phi)\right\}\rmd\phi_\Lambda\delta_\psi(d\phi_{\zd\setminus\Lambda}),
 \end{equation}
with
$$Z_{\Lambda,\psi}=\int_{\RR^\zd}\exp\left\{-\beta H_{\Lambda,\psi}(\phi)\right\}\rmd\phi_\Lambda\delta_\psi(d\phi_{\zd\setminus\Lambda}),$$
and where $\delta_\psi(d\phi_{\zd\setminus\Lambda})=\prod_{x\in\zd\setminus\Lambda}\delta_{\psi(x)}(d\phi(x))$ and determines the boundary condition. Thus, $\nu_{\Lambda,\psi}$ is characterized by the inverse temperature $\beta>0$ and the Hamiltonian $H_{\Lambda,\psi}$ on $\Lambda$, which we assume to be of gradient type:
\begin{eqnarray}
\label{eqn00}
H_{\Lambda,\psi}(\phi)=\sum_{i\in I}\sum_{x,x+e_i\in\Lambda}U(\btd_{i}\phi(x))+2\sum_{i\in I}\sum_{x\in\Lambda,x+e_i\in\partial\Lambda}U(\btd_{i}\phi(x)),
\end{eqnarray}
where the sum inside $\Lambda$ is over ordered nearest neighbours pairs $(x,x+e_i)$. We denoted by 
$$I=\{-d,-d+1,\ldots, d\}\setminus \{0\}$$
and we introduced for each $x\in\zd$ and each $i\in I$, the discrete gradient
$$\nabla_i\phi(x)=\phi(x+e_i)-\phi(x),$$
that is, the interaction depends only on the differences of neighboring heights. Note that $e_i, i=1,2,\ldots d$ denote the unit vectors and $e_{-i}=-e_i$. A model with such a Hamiltonian as defined in (\ref{eqn00}), is called a massless model
with a continuous symmetry (see \cite{FSL}). 
The potential $U\in C^2(\Bbb R)$ is a symmetric function with quadratic growth at infinity:
\begin{equation*}
\label{tag2}
U(\eta)\ge A|\eta|^2-B,\qquad \eta\in\Bbb R
\tag{A0}
\end{equation*}
for some $A>0, B\in\Bbb R$.

\subsection{General definitions and notation}

\subsubsection{$\phi$-Gibbs Measures}

For $A\subset\zd$, we shall denote by $\F_{A}$ the $\sigma$-field generated by $\{\phi(x):x\in A\}$. 
\begin{defn}  {\bf ($\phi$-Gibbs measure on $\zd$)}
\label{gibbs}
The probability measure $\nu\in P(\RR^{\zd})$ is called a Gibbs measure for the $\phi$-field with given Hamiltonian $H:=(H_{\Lambda, \psi})_{\Lambda\subset\zd, \psi \in
\RR^{\zd}}$ ($\phi$-Gibbs measure for short), if its
conditional probability of $\F_{\Lambda^c}$ satisfies the DLR equation
$$\nu(\,\cdot\,|\F_{\Lambda^c})(\psi)=\nu_{\Lambda,\psi}(\cdot),~~\nu-\mbox{a.e. }\psi,$$
for every finite $\Lambda\subset\zd$.
\end{defn}

It is known that the $\phi$-Gibbs measures exist under condition (\ref{tag2}) when the dimension $d\ge3$,
but not for $d=1,2$, where the field "delocalizes" as $\Lambda\nearrow\zd$ (see \cite{FP}). An infinite volume limit (thermodynamic limit) for $\nu_{\Lambda,\psi}$ when $\Lambda\nearrow\zd$ exists only when $d\ge3$.

\subsubsection{$\nabla\phi-$Gibbs Measures}

\textbf{Notation for the Bond Variables on $\zd$}\\[-2mm]

Let $$\bzd:=\{b=(x_b,y_b)~|~x_b,y_b\in\zd,\|x_b-y_b\|=1,b~\mbox{directed from}~x_b~\mbox{to}~y_b\};$$ 
note that each undirected bond appears twice in $\bzd$. Let 
$$\bzdl:=\bzd\cap (\Lambda\times\Lambda),~\partial\bzdl:=\{b=(x_b,y_b)~|~x_b\in\zd\setminus\Lambda,y_b\in\Lambda,\|x_b-y_b\|=1\}$$
and 
$$\bzdlb:=\{b=(x_b,y_b)\in\bzd~|~x_b\in\Lambda~\mbox{or}~y_b\in\Lambda\}.$$

For $\phi=(\phi(x))_{x\in\zd}$ and $b=(x_b,y_b)\in\bzd$, we define the \textit{height differences} $\nabla\phi(b):=\phi(y_b)-\phi(x_b)$.
The height variables $\phi=\{\phi(x);x\in\zd\}$ on $\zd$ automatically determines a field of height differences
$\nabla\phi=\{\nabla\phi(b);b\in\bzd\}$. One can therefore consider the distribution $\mu$ of $\nabla\phi$-field
under the $\phi$-Gibbs measure $\nu$. We shall call $\mu$ the $\nabla\phi$-Gibbs measure. In fact, it is possible to define
the $\nabla\phi$-Gibbs measures directly by means of the DLR equations and, in this sense, $\nabla\phi$-Gibbs measures exist for
all dimensions $d\ge1$.

A sequence of bonds $\C=\{b^{(1)},b^{(2)},\ldots,b^{(n)}\}$ is called a \textit{chain} connecting $x$ and $y$, $x,y\in\zd$, if $x_{b_1}=x,y_{b^{(i)}}=x_{b^{(i+1)}}$ for $1\le i\le n-1$ and $y_{b^{(n)}}=y$. The chain is called a \textit{closed loop} if $y_{b^{(n)}}=x_{b^{(1)}}$. A \textit{plaquette} is a closed loop $\A=\{b^{(1)},b^{(2)},b^{(3)},b^{(4)}\}$ such that $\{x_{b^{(i)}},i=1,\ldots,4\}$ consists of 
$4$ different points. 

The field $\eta=\{\eta(b)\}\in\RR^{\bzd},b\in\bzd,$ is said to satisfy \textit{the plaquette conditions} if
\begin{eqnarray}
\eta(b)=-\eta(-b)~\mbox{for all}~b\in\bzd~\mbox{and}~\sum_{b\in\A}\eta(b)=0~\mbox{for all plaquettes}~\A~\mbox{in}~\zd,
\end{eqnarray}
where $-b$ denotes the reversed bond of $b$. Let 
\begin{equation}
\chi=\{\eta\in\RR^{(\zd)^*}~\mbox{which satisfy the plaquette 
condition}\}
\end{equation}
and let $L_r^2, r>0$, be the set of all $\eta\in\RR^{\bzd}$ such that
$$|\eta|^2_r:=\sum_{b\in\bzd}|\eta(b)|^2e^{-2r\|x_b\|}<\infty.$$
We denote $\chi_r=\chi\cap L_r^2$ equipped with the norm $|\cdot|_r$. For $\phi=(\phi(x))_{x\in\zd}$ and $b\in\bzd$, we define $\eta^\phi(b):=\nabla\phi(b)$. Then $\nabla\phi=\{\nabla\phi(b)\}$ satisfies the plaquette condition. Conversely, the heights $\phi^{\eta,\phi(0)}\in\RR^\zd$ can be constructed from height differences $\eta$ and the
height variable $\phi(0)$ at $x=0$ as 
\begin{equation}
\label{19}
\phi^{\eta,\phi(0)}(x):=\sum_{b\in\C_{0,x}}\eta(b)+\phi(0),
\end{equation} 
where $\C_{0,x}$ is an arbitrary chain connecting $0$ and $x$. Note that $\phi^{\eta,\phi(0)}$ is well-defined if
$\eta=\{\eta(b)\}\in\chi$.\\[-2mm]

\textbf{Definition of $\nabla\phi$-Gibbs measures}\\[-2mm]

We next define the finite volume $\nabla\phi$-Gibbs measures. For every $\xi\in\chi$ and finite $\Lambda\subset\zd$ the space of all
possible configurations of height differences on $\bzdlb$ for given
boundary condition $\xi$ is defined as
$$\chi_{\bzdlb,\xi}=\{\eta=(\eta(b))_{b\in\bzdlb};\eta\vee\xi\in\chi\},$$
where $\eta\vee\xi\in\chi$ is determined by $(\eta\vee\xi)(b)=\eta(b)$ for $b\in\bzdlb$ and $=\xi(b)$ for
$b\not\in\bzdlb$. 
\begin{rem}\normalfont
\label{equivzd}
Note that  when $\zd\setminus\Lambda$ is connected, $\chi_{\bzdlb,\xi}$ is an affine space such that $\dim \chi_{\bzdlb,\xi}=|\Lambda|$. Indeed, fixing a point $x_0\notin\Lambda$, we consider the map $\chi_{\bzdlb,\xi}\rightarrow\RR^{\Lambda}$, such that $\eta\rightarrow\phi=\{\phi(x)\}\in\RR^{\Lambda}$, with $\phi(x)$ defined by
$$\phi(x)=\sum_{b\in C_{x_0,x}}(\eta\vee\xi)(b)$$
for a chain $C_{x_0,x}$ connecting $x_0$ and $x\in\Lambda$. This map then well-defined and an invertible linear transformation. 
\end{rem}
\begin{defn}{\bf(Finite Volume $\nabla\phi$-Gibbs measure)}
\label{finvolgrad}
The finite volume $\nabla\phi$-Gibbs measure in $\Lambda$ (or more precisely, in $\bzdl$) with given Hamiltonian $H:=(H_{\Lambda, \xi})_{\Lambda\subset\zd,\,\xi \in
\chi}$ and 
with boundary condition $\xi$ is defined by
$$
\mu_{\Lambda,\xi}(\ormd\eta)=\frac{1}{Z_{\Lambda.\xi}}\exp\left\{-\beta\sum_{b\in\bzdlb}U(\eta(b))\right\}\ormd\eta_{\Lambda,\xi}
\in P(\chi_{\bzdlb,\xi}),
$$
where $\ormd\eta_{\Lambda,\xi}$ denotes the Lebesgue measure on the affine space $\chi_{\bzdlb,\xi}$ and $Z_{\Lambda,\xi}$ is
the normalization constant.
\end{defn}
 Let $P({\chi})$ be the set of all probability measures on ${\chi}$ and let $P_2({\chi})$ be those 
$\mu\in P({\chi})$ satisfying $E^{\mu}[|\eta(b)|^2]<\infty$ for each $b\in\bzd$.
\begin{rem}\normalfont
\label{equivzd1}
For every $\xi\in\chi$ and $a\in\RR$, let $\psi=\phi^{\xi,a}$ be defined by (\ref{19}) and consider the measure $\nu_{\Lambda,\psi}$. Then $\mu_{\Lambda,\xi}$ is the image measure of $\nu_{\Lambda,\psi}$ under the map $\{\phi(x)\}_{x\in\Lambda}\rightarrow\{\eta(b):=\nabla(\phi\vee\psi)(b)\}_{b\in\bzdlb}$ and where we defined $(\phi\vee\psi)(x):=\phi(x)$ for $x\in\Lambda$ and $(\phi\vee\psi)(x):=\psi(x)$ for $x\notin\Lambda$.
Note that the image measure is determined only by $\xi$ and is independent of the choice of $a$. Let ${K}^{\psi}_{\Lambda}:\{\phi(x)\}_{x\in\zd}\rightarrow\{\eta(b)\}_{b\in\bzd}$, with $\eta(b):=\nabla(\phi\vee\psi)(b)$.
\end{rem}
\begin{defn} {\bf ($\nabla\phi$-Gibbs measure on $(\zd)^*$)} 
\label{nablaphigib}
The probability measure $\mu\in P(\chi)$ is called a Gibbs measure for the height differences with given Hamiltonian $H:=(H_{\Lambda, \xi})_{\Lambda\subset\zd, \xi \in\chi}$ ($\nabla\phi$-Gibbs measure for short),
if it satisfies the DLR equation
\begin{equation}
\label{dlr}
\mu(\,\cdot\,|\F_{\bzd\setminus\bzdlb})(\xi)=\mu_{\Lambda,\xi}(\cdot),~~\mu-\mbox{a.e. }\xi,
\end{equation}
for every finite $\Lambda\subset\zd$, where $\F_{\bzd\setminus\bzdlb}$ stands for the $\sigma$-field of $\chi$ generated by
$\{\eta(b),b\in\bzd\setminus\bzdlb\}$.
\end{defn}
\begin{rem}
Proving the DLR equation (\ref{dlr}) is equivalent to proving that for every finite $\Lambda\subset\zd$ and for all $F\in C_b(\chi)$ we have
\begin{equation}
\label{dlrgrad}
\int_\chi\mu(d\xi)\int_{\chi_{\bzdlb,\xi}} \mu_{\Lambda,\xi}(d\eta)F(\eta)=\int_\chi\mu(d\eta)F(\eta).
\end{equation}
(For a proof of this equivalence, see Remark 1.24 from \cite{Ge}).
\end{rem}
With the notations from (\ref{eqn00}) and Definition \ref{finvolgrad}, let
\begin{equation*}
{\cal G}_\beta(H):=\{\mu\in P_2(\chi):\mu~\mbox{is}~\nabla\phi-\mbox{Gibbs measure on}~\bzd~\mbox{with given Hamiltonian}~H\}.
\end{equation*}
\begin{rem}
Throughout the rest of the paper, we will use the notation $\phi,\psi$ to denote height variables and $\eta,\xi$ to denote height differences.
\end{rem}

\textbf{Shift-invariance and ergodicity}\\[-2mm]

For $x\in\zd$, we define the shift operators: $\sigma_{x}:\RR^{\zd}\rightarrow\RR^{\zd}$ for the heights by $\sigma_{x}\phi(y)=\phi(y-x)~\mbox{for}~y\in\zd~\mbox{and}~\phi\in\RR^{\zd}$, and $\sigma_{x}:\RR^{\bzd}\rightarrow\RR^{\bzd}$ for the bonds by $(\sigma_{x}\eta)(b)=\eta(b-x)$, for $b\in\bzd~\mbox{and}~\eta\in\chi$. Then \textbf{shift-invariance} and \textbf{ergodicity} for $\mu$ (with respect to $\sigma_x$ for all $x\in\zd$) is defined in the usual way (see for example page 122 in \cite{FSL}). We say that the shift-invariant $\mu\in P_2({\chi})$ has a given \textbf{tilt} $u\in\RR^d$ if $\ey_{\mu}(\eta(b))= \langle u, y_b-x_b \rangle$ for all bonds $b=(x_b,y_b)\in\bzd$.\\[-2mm]

\subsection{Results}

Our state space $\RR^{\zd}$ being unbounded, gradient interface models experience delocalization
in lower dimensions $d=1,2$, and no infinite volume Gibbs state exists in these dimensions (see \cite{FP}).
Instead of looking at the Gibbs measures of the $(\phi(x))_{x\in\zd}$,
Funaki and Spohn proposed to consider the distribution of the gradients $\left(\nabla_i\phi(x)\right)_{i\in I, x\in\zd}$ under $\nu$ (see Definition 1.5) in the \textbf{gradient Gibbs measures $\mu$},
which in view of the Hamiltonian (\ref{eqn00}), can also be given in terms of a Dobrushin-Landford-Ruelle (DLR) description. Note that infinite volume gradient Gibbs measures
exist in all dimensions, in particular for dimensions $1$ and $2$, which is one of the reasons that Funaki and Spohn introduced them. For a good background source on these models, see Funaki \cite{FSL}.

Assuming strict convexity of $U$:
\begin{equation}
\label{tag2'}
0<C_1\le U''\le C_2<\infty,
\end{equation}
Funaki and Spohn  showed in \cite{FS} the existence and uniqueness of ergodic gradient Gibbs measures for every fixed tilt $u\in\RR^d$, that is, if $\ey_{\mu}(\nabla_i\phi(x))=u_i$ for all nearest-neighbour pairs $(x,x+e_i)$ (see also \cite{Sh}).
Moreover, they also proved that the corresponding free energy, or surface tension, $\sigma(u)\in C^1(\RR^d)$
is convex in $u$; the surface tension, defined in section 7 of our paper, physically describes the macroscopic
energy of a surface with tilt $u$, i.e., a $d$-dimensional hyperplane located in $\RR^{d+1}$ with normal vector $(-u,1)\in\RR^{d+1}$. Both these results (ergodic
component and convexity of surface tension) were used in \cite{FS} for the derivation of the hydrodynamical limit of the Ginzburg-Landau model.

In fact under the strict convexity assumption (\ref{tag2'}) of $U$, much more is known for the gradient field. At large scales
it behaves much like the harmonic crystal or gradient free fields which is a Gaussian field with quadratic $U$.
In particular, Brydges and Yau \cite {BY} (in the case of small analytic perturbations of quadratic potentials), Naddaf and Spencer \cite{NS} (in the case of strictly convex potentials and tilt $u=0$) and Giacomin, Olla and Spohn \cite{gos} (in the case of strictly convex potentials and arbitrary tilt $u$) showed that the rescaled gradient field converges weakly as $\epsilon\searrow 0$ to a continuous homogeneous 
Gaussian field, that is
\begin{equation}
\label{cltconv}
S_{\epsilon}(f)=\epsilon^{d/2}\sum_{x\in\zd}\sum_{i\in I} (\nabla_i\phi(x)-u_i)f_i(\epsilon x) \rightarrow N(0,\Sigma^2_u(f))~~\mbox{as}~~\epsilon\rightarrow 0 ,\qquad
 f\in C_0^{\infty}(\RR^d;\RR^d),
\end{equation}
where the convergence takes place under ergodic $\mu$ with tilt $u$ (see Theorem 2.1 in Giacomin, Olla and Spohn \cite{gos} for an explicit expression of $\Sigma^2_u(f)$ in (\ref{cltconv}) in the case with arbitrary tilt and see Biskup and Spohn \cite{BS} for similar results in the non-convex case). This central limit theorem derived at standard scaling
$\epsilon^{d/2}$, is far from trivial since as shown in Delmotte and Deuschel \cite{DD}, the gradient field has slowly decaying,
non-absolutely summable  covariances
of the algebraic order
\begin{equation}
\label{tag3}
\left|\cov_{\nu}(\nabla_i\phi(x),\nabla_j\phi(y))\right|\sim\frac{C}{1+\|x-y\|^d}.
\end{equation}
All the above-mentioned results are proved under the essential assumption of strict convexity of the potential $U$, which assumption is necessary for the application of the Brascamp-Lieb inequality and of the Helffer-Sjostrand random walk representation (see \cite{FSL} for a detailed review of these methods and results). While
strict convexity is crucial for the proofs, one would expect some of these results to
be valid under more general circumstances, in particular also for some classes of non-convex potentials. However, so far there have been very few results
on non-convex potentials. This is where the focus of this paper comes in,
which is to extend the results known for strictly convex potentials to some classes of non-convex potentials.

We will briefly summarize next
the state of affairs regarding results for non-convex potentials, in the different regimes at inverse temperature $\beta$. At low temperature (i.e. large $\beta$) using the renormalization
group techniques developed by Brydges \cite{david}, Adams et al. \cite {AKM} show in on-going work for a class of non-convex potentials, the strict convexity of the surface tension for small tilt $u$. At moderate temperature ($\beta=1$), Biskup and Koteck\'y \cite{BK} give an example 
of a non-convex potential $U$ for which uniqueness of the ergodic gradient Gibbs measures $\mu$ fails. The potential $U$ can be described as the mixture of two Gaussians with two different variances. For this particular case of $U$, \cite{BK} prove co-existence of two ergodic gradient Gibbs measures at tilt $u=0$ (see also Figure 4 and example 3.2 (a) below). See also the work of
Fr\"ohlich and Spencer (\cite{CMP81}, \cite{JSP81}) in relation to the Coulomb gas, and
the theory based on the infrared-bound (e.g. Fr\"ohlich, Simon and
Spencer \cite{CMP76}). For high temperature (i.e. small $\beta$), we have proved in a previous paper with S. Mueller \cite{CDM} strict convexity of the surface tension in a regime similar to (\ref{tag5}) below. Our potentials are of the form
$$U(\nabla_i\phi(x))=V(\nabla_i\phi(x))+g(\nabla_i\phi(x))$$
where $V,g\in C^2(\RR)$ are such that 
\begin{equation*}
C_1\le V''\le C_2,~0<C_1<C_2~~~~\mbox{and}~~~~-C_0\le g''\le 0,~\mbox{with}~C_0>C_2.
\end{equation*}
Specifically, we assumed in \cite{CDM} that
$$\frac{4}{\pi}(12d\bar{C})^{1/2}\sqrt{\beta C_1}\frac{1}{C_1}||g''||_{L^1(\RR)}\le\frac{1}{2},~\mbox{where}~\bar{C}=\max\left(\frac{C_0}{C_1},\frac{C_2}{C_1}-1,1\right).$$
The method used in \cite{CDM}, based on two scale decomposition of the free field,
gives less sharp estimates for the temperature than our current paper as the estimates also depend on $C_0$. However, at this point it is not clear whether the method introduced in \cite{CDM} could yield any other result of interest than the strict convexity of the surface tension.

The aim of our current paper is to use an alternative technique from the one we used in \cite{CDM} and relax the strict convexity assumption (\ref{tag2'}) to obtain much more than just strict convexity of the surface tension; more precisely, we also prove uniqueness of the ergodic component
at every tilt $u\in\RR^d$, central limit theorem of form as given in (\ref{cltconv}) and decay of covariances as in (\ref{tag3}). As stated above, the hydrodynamical limit for the corresponding Ginzburg-Landau model should then essentially
follow from our results. Our main results are proven under the assumption that 
\begin{equation*}
\label{vc}
C_1\le V''\le C_2,~0<C_1<C_2~~~~\mbox{and}~~~~-\infty< g''\le 0
\tag{A1}
\end{equation*}
and that the inverse temperature $\beta$ is sufficiently small, that is if
\begin{equation*}
\label{tag5}
\beta^{\frac{1}{2q}}||g''||_{L^q(\RR)}<\frac{(C_1)^{\frac{3}{2}}}{2C_2^{\frac{q+1}{2q}}\left(2d\right)^{\frac{1}{2q}}},~\mbox{for some}~q\ge 1,
\tag{A2}
\end{equation*}
or if
\begin{equation*}
\label{g'l2}
\beta^{\frac{3}{4}}||g'||_{L^2(\RR)}\le \frac{(C_1)^{\frac{3}{2}}}{2(C_2)^{\frac{5}{4}}(2d)^{\frac{3}{4}}}.
\tag{A3}
\end{equation*}
The condition (\ref{vc}) with $g''\le 0$ may look a bit artificial, but as we elaborate in Remark \ref{compsup} in section 3 below, any perturbation $g\in C^2$ with compact support can be substituted for the $g''\le 0$ assumption in (\ref{vc}). Note that in contrast to the condition in our previous paper \cite{CDM}, $||g''||_{L^{\infty}(\RR)}$ can be arbitrarily large as long as $||g''||_{L^q(\RR)}$ is small. Note also that using an obvious rescaling argument (see Remark \ref{scal}), we can always reduce our
assumption (\ref{vc}) to the
case $\beta=C_1=1$; then (\ref{tag5}), respectively (\ref{g'l2}), states that our condition is satisfied
whenever the perturbation $g''$
is small in the $L^q(\RR)$, respectively $g'$ is small in the $L^2(\RR)$ sense.

Our main result is the following 
\begin{thm} [Uniqueness of an ergodic $\mu_u$]
\label{ergall} 
Let $U=V+g$, where $U$ satisfy (\ref{tag2}) and $V$ and $g$ satisfy (\ref{vc}) and (\ref{tag5}) or (\ref{vc}) and (\ref{g'l2}). Then for every $u\in\RR^d$, there exists at most one ergodic, shift-invariant $\mu_u\in {\cal G}_\beta(H)$ with a given tilt $u\in\RR^d$.
\end{thm} 
Let $F\in C^1_b(\chi_r)$, where $C^1_b(\chi_r)$ denotes the set of differentiable functions depending on finitely many 
coordinates with bounded derivatives and where $\chi_r$ was defined in subsection 1.2.2. For $\eta,\eta'\in\chi$, let
$$\lim_{\epsilon\rightarrow 0}\frac{F(\eta+\epsilon\eta')-F(\eta)}{\epsilon}=\langle DF(\eta),\eta'\rangle=\sum_{b\in\bzd}\alpha(b)\eta'(b).$$
We denote by
\begin{equation}
\label{partialbond}
\partial_bF(\eta):=\alpha(b)~\mbox{and}~||\partial_bF||_{\infty}=\sup_{\eta\in\chi}|\partial_bF(\eta)|.
\end{equation}
Another result we prove for our class of non-convex potentials is
\begin{thm} [Decay of Covariances]
\label{cov} 
Let $u\in\RR^d$. Assume $U=V+g$, where $U$ satisfies (\ref{tag2}) and $V$ and $g$ satisfy (\ref{vc}) and  (\ref{tag5}) or (\ref{vc}) and (\ref{g'l2}). Let $F,G\in C^1_b(\chi_r)$. Then there exists $C>0$ such that
\begin{equation}
\label{covzd}
|\cov_{\mu_u}(F(\eta),G(\eta))|\le C\sum_{b,b'\in\bzd}\frac{||\partial_b F||_{\infty}||\partial_{b'}G||_{\infty}}{1+\|x_b-x_{b'}\|^d},
\end{equation}
where $b=(x_b,y_b)$ and $b'=(x_{b'},y_{b'})$.
\end{thm}
We also prove
\begin{thm} [Central Limit Theorem]
\label{clt} 
Let $u\in\RR^d$. Assume $U=V+g$, where $U$ satisfies (\ref{tag2}) and $V$ and $g$ satisfy (\ref{vc}) and  (\ref{tag5}) or (\ref{vc}) and (\ref{g'l2}).
Set
$$S_{\epsilon}(f)=\epsilon^{d/2}\sum_{x\in\zd}\sum_{i\in I} (\nabla_i\phi(x)-u_i)f_i(\epsilon x),$$
where $f\in C_0^{\infty}(\RR^d;\RR^d)$. Then
$$S_{\epsilon}(f)\Rightarrow N(0,\Sigma^2_u(f))~~\mbox{as}~~\epsilon\rightarrow 0,$$
where $\Sigma^2_u(f)$ can be identified explicitly as in Theorem 2.1 in \cite{gos}, $\Sigma^2_u(f)\neq 0$ for $f\neq 0$, and $\Rightarrow$ signifies convergence in distribution.
 \end{thm}
Moreover, we extend in Theorem \ref{evconv1} the results of strict convexity of the surface tension from \cite{FS} and \cite{DGI} to the family of non-convex potentials satisfying (\ref{tag2}), (\ref{vc}) and (\ref{tag5}).
 
Even though our results are obtained for the high
temperature case, previously only our results in \cite{CDM} were known for the non-convex case. Also, the proofs of this paper require some crucial observations not made before. Moreover, in our main result Theorem \ref{ergall}, we prove uniqueness of ergodic gradient Gibbs measures $\mu$ with a given arbitrary tilt $u\in\RR^d$ for the class of non-convex potentials satisfying (\ref{tag2}), (\ref{vc}) and (\ref{tag5}). To the best of our knowledge, this is the first result where uniqueness of ergodic gradient Gibbs measures $\mu$ is proved for a class of non-convex potentials $U$. For potentials that are mixtures of Gaussians as considered in Biskup and Koteck\'y \cite{BK}, they prove non-uniqueness of ergodic gradient Gibbs measures for tilt $u=0$ in the $\beta=1$ regime. For the same example, we prove uniqueness of ergodic gradient Gibbs measures for given arbitrary tilt $u$ in the high temperature regime. Therefore, our result also highlights the existence of phase transition for these models in different temperature regimes.

The basic idea relies on a one-step coarse graining procedure, in which we consider the
marginal distribution of the gradient field restricted to the even sites,
which is also a gradient Gibbs field.
The corresponding Hamiltonian, although no longer a two-body Hamiltonian, is
then obtained via integrating out the field at
the odds sites. We can integrate out the field $\phi$ at all odd
sites, using the fact that they are independent conditional on the field $\phi$ at
even sites, which is a consequence of the bi-partiteness of the graph $\Z^d$ with nearest-neighbor bonds. The crucial step, which is similar to the idea of our previous paper \cite{CDM}, is that
strict convexity can be gained via integration
at sufficiently high temperature (see also Brascamp et al. \cite{BLL} for previous use of the even/odd representation). The essential observation is that we can
formulate a condition for this multi-body potential, which we call the {\it random walk representation condition}, which allows us to obtain a strictly convex Hamiltonian, and implies the random walk representation, permitting us to apply the techniques of Helffer and Sj{\"o}strand \cite{HS} or Deuschel \cite{JD1}. The random walk representation condition, and implicitly the strict convexity of the new Hamiltonian, can be verified under our assumptions as in (\ref{tag2}), (\ref{vc}) and (\ref{tag5}). Note that the method in \cite{CDM} is more general and could be applied to non-bipartite graphs. 

A natural
question to ask is whether we can iterate the coarse graining procedure in our current paper and find a scheme which could possibly lower the temperature towards the critical $\beta_c$, which marks the transition from a unique gradient Gibbs measure $\mu$ (as proved in Theorem \ref{ergall} in our paper for arbitrary tilt $u$) to multiple gradient Gibbs measures $\mu$ (as proved in \cite{BK} for tilt $u=0$). Note that iterating the coarse graining scheme is an interesting open problem. One of the
main difficulties is that, after iteration, the bond structure on the even
sites of $\Z^d$ changes, and we no longer have a bi-partite graph. Currently, we could use our method as detailed in sections 2 and 3, to keep integrating out lattice points so that the new Hamiltonian at each step, always of gradient type, can be separated into a strictly convex part and a non-convex perturbation; however, at this point, our technique for estimation of covariances as given in section 3, is not robust enough to allow us to keep coarse graining the lattice points for more than a finite number of steps, before we stop being able to improve the assumptions on our initial perturbation $g$.

The rest of the paper is organized as follows: In section 2 we present the odd/even characterization of the gradient field. In section 3 we give the formulation of the random walk representation condition, which is verified in Theorem \ref{rw} under conditions (\ref{tag2}), (\ref{vc}) and (\ref{tag5}). Section 3 also presents a few examples, in particular we show that our criteria gets close
to the Biskup-Koteck\'y phase co-existence regime, both for the case of the zero and the non-zero tilt $u$ (see example 3.2 (a)). In section 4 we prove Theorem \ref{ergall}, our main result on uniqueness of ergodic gradient Gibbs measure with given tilt $u$, which is based on adaptations of \cite{FS}, assuming the random walk representation condition. Section 5 deals with the decay of covariances and the proof is based on the random walk representation for the field at the even sites which allows us to use the result of \cite{DD}. Section 6 shows the central limit theorem, here again
we focus on the field at even sites and apply the random walk representation idea of \cite{gos}. Section 7 proves the strict convexity of the surface tension, or free energy, which follows from the
convexity of the Hamiltonian for the gradient field restricted to the even sites. Finally, the appendix provides explicit computations for our one-step coarse graining procedure in the special case of potentials considered by \cite{BK} (see also example 3.2 (a)).

\section{Even/Odd Representation}

There are two key results in this section.  The first one is Lemma \ref{muodd}, where we are restricting the height differences
to the even sites, which induces a $\nabla\phi$ measure on the even lattice with
a different bond structure.  The second main result of this section is Lemma \ref{led}, where we give a formula for the conditional
of a $\nabla\phi$-Gibbs measure on the height differences between even sites. These two results will be essential for the proof for one of our main results, that is for the proof of the uniqueness of ergodic component of Theorem \ref{ergall}.

In Subsection 2.1 we introduce the notation for the bond variables on the even subset of $\zd$, in Subsection 2.2 we define the $\phi$-Gibbs measure and the $\nabla\phi$-Gibbs measure corresponding to the even subset of $\zd$ and in Subsection 2.3 we present the relationship between the $\nabla\phi$-Gibbs measures for the bonds on $\zd$ and the $\nabla\phi$ for the bonds on even subset of $\zd$, when their corresponding finite volume $\phi$-Gibbs measures are related by restriction. 

\subsection{Notation for the Bond Variables on the Even Subset of $\zd$}
As $\zd$ is a bipartite graph, we will label the vertices of $\zd$ as \textbf{even} and \textbf{odd} vertices, such that every \textbf{even} vertex has only \textbf{odd} nearest neighbor vertices and vice-versa. 

Let 
$$\ed:=\{a=(a_1,a_2,\ldots,a_{d})\in\zd~|~\sum_{i=1}^d a_i=2p,p\in\Z\}$$ 
and 
$$\od:=\{a=(a_1,a_2,\ldots,a_{d})\in\zd~|~\sum_{i=1}^d a_i=2p+1,p\in\Z\}.$$ Let $\Lambda_\even\subset\ed$ finite. We will next define the bonds in $\ed$ in a similar fashion to the definitions for bonds on $\zd$. Let 
$$\bed:=\{b=(x_b,y_b)~|~x_b,y_b\in\ed,\|x_b-y_b\|=2,b~\mbox{directed from}~x_b~\mbox{to}~y_b\},$$ 
$$\bedl:=\bed\cap (\Lambda_\even\times\Lambda_\even),~~\overline{\bedl}:=\{b=(x_b,y_b)\in\bed~|~x_b\in\Lambda_\even~\mbox{or}~y_b\in\Lambda_\even\},$$
$$\partial\bedl:=\{b=(x_b,y_b)~|~x_b\in\ed\setminus\Lambda_\even,y_b\in\Lambda_\even,\|x_b-y_b\|=2\}$$ and 
$$\partial\Lambda_\even:= \left\{y\in\ed\setminus\Lambda_\even~|~,\|y-x\|=2 ~\mbox{for some}~x\in\Lambda_\even\right\}.$$ 
Note that throughout the rest of the paper, we will refer to the bonds on $\bed$ as the \textbf{even bonds}.  

\begin{figure}
\centering
\resizebox{!}{3.5cm}{%
\includegraphics{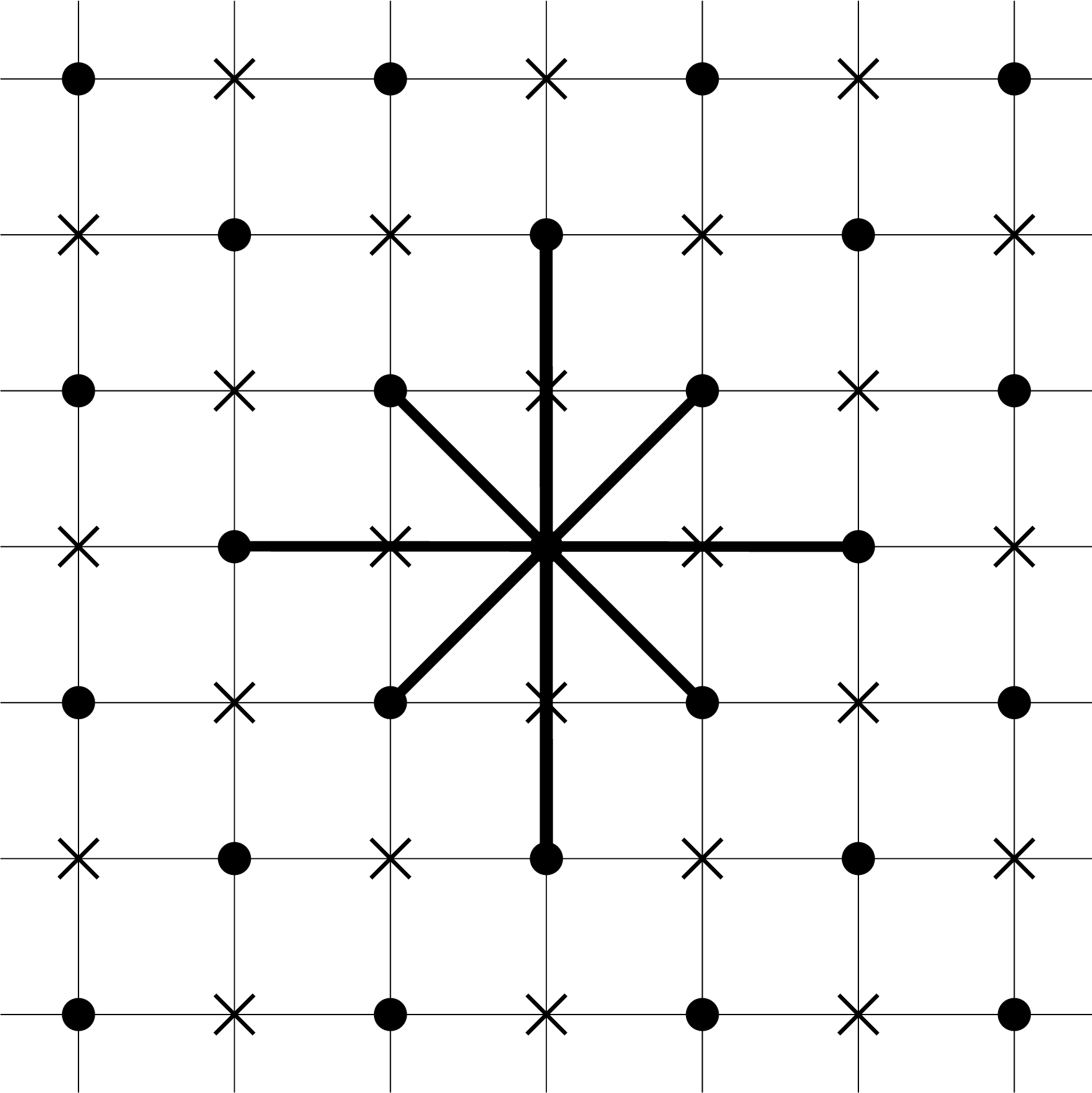}}
\caption{
The bonds of $0$ in $Z^2_{\even}$}
\end{figure}

An \textbf{even} \textit{plaquette} is a closed loop $\A_\even=\{b^{(1)},b^{(2)},\ldots,b^{(n)}\}$, where $b^{(i)}\in\bed$, $n\in\{3,4\}$, such that $\{x_{b^{(i)}},i=1,\ldots,n\}$ consists of 
$n$ different points in $\ed$. 
The field $\eta=\{\eta(b)\}\in\RR^{\bed}$ is said to satisfy the \textbf{even} \textit{plaquette condition} if
\begin{equation}
\eta(b)=-\eta(-b)~\mbox{for all}~b\in\bed~\mbox{and}~\sum_{b\in{\A_\even}}\eta(b)=0~\mbox{for all \textbf{even} plaquettes in}~\ed.
\end{equation}
Let $\chi_\even$ be the set of all $\eta\in\RR^{\bed}$ which satisfy the even plaquette 
condition.
For each $b=(x_b,y_b)\in \bed$ we define the \textit{even height differences} $\eta_\even(b):=\nabla_\even\phi(b)=\phi(y_b)-\phi(x_b)$.
The heights $\phi^{\eta_\even,\phi(0)}$ can be constructed from the height differences $\eta_\even$ and the height variable $\phi(0)$ at $x=0$ as 
\begin{equation}
\label{19ev}
\phi^{\eta_\even,\phi(0)}(x):=\sum_{b\in C^\even_{0,x}}\eta_\even(b)+\phi(0),
\end{equation}
where $x\in\ed$ and $C^\even_{0,x}$ is an arbitrary path in $\ed$ connecting $0$ and $x$. Note that $\phi^{\eta,\phi(0)}(x)$ is well-defined if $\eta_\even=\{\eta_\even(b)\}\in{\chi}_\even$. We also define $\chi_{\even,r}$ similarly as we define $\chi_r$. As on $\zd$, let $P({\chi}_\even)$ be the set of all probability measures on ${\chi}_\even$ and let 
$P_2({\chi}_\even)$ be those $\mu\in P({\chi}_\even)$ satisfying $E^{\mu}[|\eta_\even(b)|^2]<\infty$ for each $b\in\bed$. We denote $\chi_{\even,r}=\chi\cap L_r^2$ equipped with the norm $|\cdot|_r$.
\begin{rem}
\label{chitochie}
Let $\eta\in\chi$. Using the plaquette condition property of $\eta$, we will define $\eta_\even$, the induced bond variables on the even lattice, from $\eta$ thus: if $b_1=(x,x+e_i)$, $b_2=(x+e_j,x)$ and $b_\even=(x+e_j,x+e_i)$, we define $\eta_\even(b_\even)=\eta(b_1)+\eta(b_2)$. Note that $\eta_\even\in\chi_\even$.
\end{rem}
\begin{rem}
Throughout the rest of the paper, we will use the notation $\phi_\even,\psi_\even$
either for a stand alone configuration on the even vertices, or in relation to the restriction of $\phi$ to the even vertices. $\eta_\even, \xi_\even$ will denote configurations on the even bonds. Similarly, $\Lambda_{\even}$ will either be a stand alone subset of $\ed$ or will be used in relation to the restriction of a set
$\Lambda\subset\zd$ to $\ed$.
For $\Lambda\subset\zd$, we will denote $\odl:=\od\cap\Lambda$.
\end{rem}

\subsection{Definition of $\nabla\phi$-Gibbs measure on $\bed$}
For every $\xi_\even\in\chi_\even$ and finite $\Lambda_\even\subset\ed$, the space of all
possible configurations of height differences on $\overline{\bedl}$ for given
boundary condition $\xi_\even$ is defined as
$$\chi_{\overline{\bedl},\xi_\even}=\{\eta_\even=(\eta_\even(b))_{b\in\overline{\bedl}},\eta_\even\vee\xi_\even\in\chi_\even\},$$
where $\eta_\even\vee\xi_\even\in\chi_\even$ is determined by $(\eta_\even\vee\xi_\even)(b)=\eta_\even(b)$ for $b\in\overline{\bedl}$ and $=\xi_\even(b)$ for $b\not\in\overline{\bedl}$. 

The $\phi$-\textit{Gibbs measure} $\nu^\even$ on $\ed$ and the $\nabla\phi$-\textit{Gibbs measure} $\mu^\even$ on $\bed$ with given Hamiltonian $H^\even$ can be defined similarly to the $\phi$-Gibbs measure and the $\nabla\phi$-Gibbs measure in Subsections 2.1 and  2.2.2. They are basically a $\phi$-Gibbs and $\nabla\phi$-Gibbs measure on a different graph, with vertex and edge
sets $(\ed,\bed)$. They are defined via the corresponding Hamiltonian $H^\even_{\Lambda_\even,\xi_{\even}}$, assumed of \textit{even} gradient type, via the \textit{finite volume Gibbs measure} $\nu^\even_{\Lambda_\even,\psi_\even}$ on $\ed$ and the \textit{finite volume $\nabla$-Gibbs measure} 
$\mu^\even_{\Lambda_\even,\psi_\even} $ on $\bed$.

Let 
$$H^\even:=(H^\even_{\Lambda_\even, \xi_\even})_{\Lambda_\even\subset\ed,\xi_\even\in\chi_\even}$$
and let
$${\cal G_\even}(H^{\even}):=\{\mu_\even\in P_2(\chi_\even):\mu^\even~\mbox{is}~\nabla\phi-\mbox{Gibbs measure on}~\bed~\mbox{with given Hamiltonian}~H^\even\}.$$
\begin{rem}\normalfont
\label{equived}
Similar to Remark \ref{equivzd}, when $\ed\setminus\Lambda_\even$ is 
connected, $\chi_{\overline{\bedl},\xi_\even}$ is an affine space such that $\dim \chi_{\overline{\bedl},\xi_\even}=|\Lambda_\even|$. Fixing a point $x_0\notin\Lambda_\even$, we consider the map ${J}^{\even,\xi}_{\Lambda_\even}:\chi_\even\rightarrow\RR^{\ed}$, such that $\eta_\even\rightarrow\{\phi_\even(x)\}$, with 
$$\phi(x):=\sum_{b\in C^{\even}_{x_0,x}}(\eta_\even\vee\xi_\even)(b),\,\,\,x\in\Lambda_\even$$ 
for a chain $C^{\even}_{x_0,x}$ connecting $x_0$ and $x$ and for fixed $\phi(x_0)$,
$$\phi(x):=\psi^{\xi_{\rm ev},\phi(x_0)}(x)=\sum_{b\in\C_{x_0,x}}\xi_\even(b)+\phi(x_0),\,\,\,x\notin\Lambda_\even.$$
\end{rem}
\begin{rem}\normalfont
\label{equived1}
For every $\xi_{\rm ev}\in\chi_{\even}$ and $a\in\RR$, let $\psi_{\even}=\phi^{\xi_{\even},a}$ be defined by (\ref{19ev}) and consider the measure 
$\nu_{\Lambda_{\rm ev},\psi_{\rm ev}}$. Then $\mu_{\Lambda_{\rm ev},\psi_{\rm ev}} $ is the image measure of $\nu_{\Lambda_{\rm ev},\psi_{\rm ev}}$ under the map $\{\phi(x)\}_{x\in\Lambda_{\rm ev}}\rightarrow\{\eta_{\rm ev}(b):=\nabla(\phi_{\rm ev}\vee\psi_{\rm ev})(b)\}_{b\in\overline{\bedl}}$.
Note that the image measure is determined only by $\xi_{\rm ev}$ and is independent of the choice of $a$.
\end{rem}

\subsection{Induced $\nabla\phi$-Gibbs measure on $(\ed)^*$} 

Throughout this section, we will make the following notation conventions. For $\phi, \psi\in\RR^{\zd}$, we define $\phi_\even:=(\phi(x))_{x\in\ed}, \psi_\even:=(\psi(x))_{x\in\ed}$. For $\eta,\xi\in\chi$, we define $\eta_\even$ and $\xi_\even$ according to Remark \ref{chitochie}.
\begin{defn}
\label{ascolambda}
Let $\Lambda_\even$ be a finite set in $\ed$. We construct a finite set $\Lambda\subset\zd$ associated to $\Lambda_\even$ as follows: if $x\in\Lambda_\even$, then $x\in\Lambda$ and $x+e_i\in\Lambda$ for all $i\in I=\{-n,-n+1,\ldots,n\}\setminus\{0\}$. Note that by definition, $\partial\Lambda=\partial\Lambda_\even$, where the boundary operations are performed in the graphs $(\zd, (\zd)^*)$ and $(\ed,\bed)$, respectively. (see Figures 2 and 3). 
\end{defn}
\begin{figure}[b]
~\hfill
\begin{minipage}[t]{5cm}
\resizebox{!}{3.5cm}{%
\includegraphics{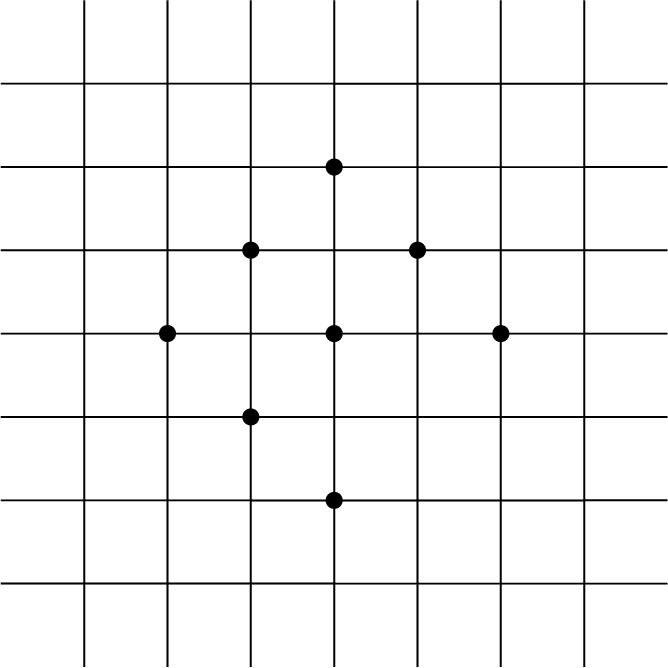}}
\caption{The graph of $\Lambda_\even$}
\end{minipage}
\hfill
\begin{minipage}[t]{5cm}
\resizebox{!}{3.5cm}{%
\includegraphics{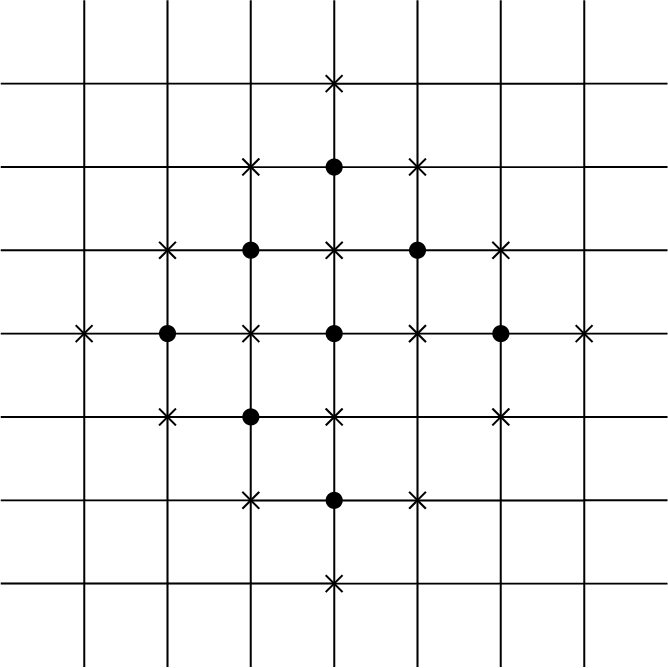}}
\caption{The graph of $\Lambda$ associated to $\Lambda_\even$}
\end{minipage}
\hfill~
\end{figure}
\begin{lemma} [Induced finite volume $\phi$-Gibbs measure on $\ed$]
\label{nuodd} 
Let $\Lambda_\even\subset\ed$ and let $\Lambda$ be the associated set in $\zd$, as defined in Definition \ref{ascolambda}. Let $\nu_{\Lambda,\psi}$ be the finite volume Gibbs measure on $\Lambda$ with boundary condition $\psi$ and with Hamiltonian $H_{\Lambda,\psi}$ defined as in (\ref{eqn00}). We define the induced finite volume Gibbs measure on $\ed$ as $\nu^\even_{\Lambda_\even,{\psi}_\even}:=\nu_{\Lambda,\psi}|_{\F(\ed)}$. Then $\nu^\even_{\Lambda_\even,{\psi}_\even}$ has Hamiltonian $H^\even_{\Lambda_\even, \psi_\even}$, where
\begin{equation}
\label{eqnW}
\begin{array}{cc}
&H^\even_{\Lambda_\even,\psi_\even}(\phi_\even):=\sum_{x\in\odl}F_x((\phi(x+e_i))_{i\in I}),\\
\mbox{\normalfont with}&\\
&F_x((\phi(x+e_i))_{i\in I})=-\log\int_{\RR} e^{-2\beta\sum_{i\in I} U(\btd_{i}\phi(x))}\rmd\phi(x).
\end{array}
\end{equation}
\end{lemma}
\begin{rem}\normalfont
\label{rgrad}
Note that for any constant $C\in\RR$, by using the change of variables $ \phi(x)\rightarrow\phi(x)+C $ in the integral formula for $F_x((\phi(x+e_i))_{i\in I})$ in (\ref{eqnW}), we have
$$F_x((\phi(x+e_i))_{i\in I})=F_x((\phi(x+e_i)+C)_{i\in I}).$$
In particular, this means that for any fixed $k\in I$
\begin{equation}
\label{weiej}
F_x((\phi(x+e_i))_{i\in I})=F_x((\phi(x+e_i)-\phi(x+e_k))_{i\in I}).
\end{equation}
Therefore we are still dealing with a \textit{gradient system}. However, it is in
general no longer a two-body gradient system. $F_x((\phi(x+e_i))_{i\in I})$, and consequently $H^{\even}_{\Lambda_\even,\psi_\even}$, are functions of the \textit{even gradients} by (\ref{weiej}) and (\ref{eqnW}).
\end{rem}
\begin{rem}\normalfont
\label{constrevham}
We formulate next more explicitly the
dependence of $F_x$ and
$H^\even_{\Lambda_\even, \psi_\even}$ on the even gradients. Let $k\in I$ be arbitrarily fixed. For any $x\in\zd$, let
$${\mathcal B}(x,k)=\{(x+e_k,x+e_i)\}_{i\in I}.$$
For all $\Lambda_\even\subset\ed$, take the set $\Lambda$ associated to $\Lambda_\even$, as defined in Definition \ref{ascolambda}.
We define here
$H^\even:=(H^\even_{\Lambda_\even, \xi_\even})_{\Lambda_\even\subset\ed,\xi_\even\in\chi_\even}$ as follows
\begin{equation}
\label{hameven}
H^\even_{\Lambda_\even, \xi_\even}(\eta)=\sum_{x\in\odl}F_x\left((\eta_\even(b))_{b\in{\mathcal B}(x,k)}\right).
\end{equation}
Note that, via Remark \ref{equived}, one can easily obtain the equivalence between the corresponding finite volume $\phi$-Gibbs and $\nabla\phi$-Gibbs measures.
\end{rem}
\begin{rem}\normalfont
By definition, $F_x((\phi(x+e_i))_{i\in I})$ only depend on sites within distance $2$ of $x$. Note that the new Hamiltonian $H_{\Lambda_\even,\psi_\even}$ depends on $\beta$ through the functions $F_x((\phi(x+e_i))_{i\in I})$. 
\end{rem}
\textbf{Proof of Lemma \ref{nuodd}}
The idea of this proof is just integrating out the height
variables on the odd sites, conditioned on the even sites. The Gibbs
property and specific
graph structure imply that the odd height variables are independent
conditional on the even sites.

Set 
\begin{equation}
\label{hxeqn}
H_x(\phi)=\sum_{i\in I} U(\nabla_i\phi(x)).
\end{equation}
Let $\Lambda_\even$ be a finite set in $\ed$ and let $\Lambda\in\zd$ be the associated set as defined in Definition \ref{ascolambda}.
Note now that due to the symmetry of the potential $U$, to the specific boundary conditions on $\Lambda$ and by (\ref{eqn00}), we have
\begin{equation}
\label{hex1}
H_{\Lambda,\psi}(\phi)=\sum_{x\in\bar\Lambda} H_x(\phi)=2\sum_{x\in\odl} H_x(\phi).
\end{equation}
Let $A\in {{\cal F}}_\ed\subset {\cal F}_\zd$, $\rmd\phi_{\Lambda_\even}=\prod_{x\in {\Lambda}_\even}\rmd\phi(x)$ and  $\rmd\phi_{\odl}=\prod_{x\in\odl}\rmd\phi(x)$. Recall that $\bar\Lambda=\Lambda\cup\partial\Lambda$ and take ${\overline{\Lambda}}_\even=\bar{\Lambda}\cap\ed$ and ${\overline{\Lambda}}_\odd=\bar{\Lambda}\cap\od$.
Then, by integrating out the odd height variables conditional on the even height variables, due 
to the Gibbs property of $\nu_{\Lambda, \psi}$ (see Definition \ref{gibbs}) and to the fact that $\partial\Lambda=\partial\Lambda_\even$, we have for every $\psi\in\RR^{\zd}$
\begin{eqnarray}
\label{gp}
\lefteqn{\nu_{\Lambda,\psi}(A)}\nonumber\\
&=&\frac{1}{Z_{\Lambda,\psi}}\int_{\Bbb R^{\bar{\Lambda}}}1_{A}(\phi) e^{-\beta
H_{\Lambda,\psi}(\phi)}d\phi_{\Lambda}\delta_\psi(d\phi_{\zd\setminus\Lambda})\nonumber\\
&{ \mbox{by (\ref{hex1})} \atop = }&\frac{1}{Z_{\Lambda,\psi}}\int_{\Bbb R^{\bar{\Lambda}}}1_{A}(\phi) e^{-2\beta\sum_{x\in\odl}
H_x(\phi)}\rmd\phi_{\odl}\rmd\phi_{\Lambda_\even}\delta_\psi(d\phi_{\zd\setminus\Lambda})\nonumber\\
&=&\frac{1}{Z_{\Lambda,\psi}}\int_{\RR^{{\overline{\Lambda}}_\even}}\int_{\RR^{{\overline{\Lambda}_\odd}}}1_{A}(\phi) e^{-2\beta\sum_{x\in\odl}
H_x(\phi)}\rmd\phi_{\odl}\rmd\phi_{\Lambda_\even} \delta_\psi(d\phi_{\zd\setminus\Lambda})\nonumber\\
&{ \mbox{as $A\in {{\cal F}}_\ed$ } \atop = }&\frac{1}{Z_{\Lambda,\psi}}\int_{\RR^{{\overline{\Lambda}}_\even}}1_{A}(\phi)\left(\int_{\RR^{{\overline{\Lambda}_\odd}}} e^{-2\beta\sum_{x\in\odl}
H_x(\phi)}\rmd\phi_{\odl} \right)\rmd\phi_{\Lambda_\even} \delta_\psi(d\phi_{\zd\setminus\Lambda})\nonumber\\
&=&\frac{1}{Z_{\Lambda,\psi}}\int_{\RR^{{\overline{\Lambda}}_\even}}1_{A}(\phi)\left(\int_{\RR^{{\overline{\Lambda}_\odd}}} \prod_{x\in\odl}e^{-2\beta
H_x(\phi)}\prod_{x\in\odl}\rmd\phi(x) \right)\rmd\phi_{\Lambda_\even} \delta_\psi(d\phi_{\zd\setminus\Lambda})\nonumber\\
&{ \mbox{as $\partial\Lambda=\partial\Lambda_\even$} \atop = }&\frac{1}{Z_{\Lambda,\psi}}\int_{\RR^{{\overline{\Lambda}}_\even}}1_{A}(\phi)\prod_{x\in\odl}\left(\int_{\Bbb R} e^{-2\beta
H_x(\phi)}\rmd\phi(x)\right)\rmd\phi_{\Lambda_\even} \delta_\psi(d\phi_{\ed\setminus\Lambda_\even})\nonumber\\
&{ \mbox{by (\ref{eqnW})} \atop = }&\frac{1}{Z_{\Lambda,\psi}}\int_{A} e^{-\sum_{x\in\odl}F_x((\phi(x+e_i))_{i\in I})}\rmd\phi_{\Lambda_\even}\delta_\psi(d\phi_{\ed\setminus\Lambda_\even})=\nu^\even_{\Lambda_\even,\psi_\even}(A),
\end{eqnarray}
where for the last equality we used that $Z_{\Lambda,\psi}=Z_{\Lambda_\even,\psi_\even}$, which is due to the fact that $\partial\Lambda=\partial\Lambda_\even$.
\endpf
\begin{lemma}
\label{muodd} {\bf(Induced $\nabla\phi$-Gibbs measure on $\bed$)}
Let $\mu\in {\cal G}_\beta(H)$. We define the induced $\nabla\phi$-Gibbs measure on $\bed$ as $\mu^\even:=\mu|_{\F\bed}$. Then $\mu^\even\in {\cal G}_\even(H^\even)$, where $H^\even_{\Lambda_\even,\xi_\even} $ is defined as in Remark \ref{constrevham}.
\end{lemma}
\pf
\Let
$${{\cal F}}_\bzd:=\sigma\left(\eta(b),b\in\bzd\right)~\mbox{and}~{{\cal F}}_\bed:=\sigma\left(\eta_\even(b),b\in\bed\right).$$
To prove the statement of the theorem, we need to prove that for all $A\in {{\cal F}}_{\bed}$, $\mu^{\even}$ satisfies
$$\mu^\even(A|{{\cal F}}_{(\ed)^*\setminus {\overline{(\Lambda_\even)}}^*})(\xi_\even)=\mu^\even_{\Lambda_\even,\xi_\even}(A).$$
In order to prove the above equality, we will first show that for all $A\in {{\cal F}}_{\bed}$ and for any $\Lambda_\even$ finite set in $\ed$ with associated set $\Lambda\subset\zd$ as defined in Definition \ref{ascolambda}, we have
\begin{equation}
\label{equivfin}
\mu_{\Lambda,\xi}(A)=\mu^\even_{\Lambda_\even,\xi_\even}(A).
\end{equation}
Then using ${\cal F}_{(\ed)^*\setminus {\overline{(\Lambda^{\cal E })}}^*}\subset {\cal F}_{(\zd)^*\setminus {\overline{(\Lambda)}}^*}$, the definition of the $\nabla\phi$-Gibbs measure and (\ref{equivfin}), we have
$$\mu(A|{{\cal F}}_{(\ed)^*\setminus {\overline{(\Lambda_\even)}}^*})(\xi)=
\ey_{\mu}\left(\ey_{\mu}\left(1_A|{{\cal F}}_{(\zd)^*\setminus {\overline{(\Lambda)}}^*}\right)|{{\cal F}}_{(\ed)^*\setminus {\overline{(\Lambda_\even)}}^*}\right)(\xi)=
\mu^\even_{\Lambda_\even,\xi_\even}(A).$$
The key point in the above equation is that when we condition further, we get $\mu_{\Lambda,
\xi'}$ where $\xi'$ is random and
being integrated over, and $\xi'$ all have $\xi_\even$ as its restriction on the
evens, and for all such $\xi'$, by (\ref{equivfin})
$\mu_{\Lambda, \xi'}$ all equal $\mu^\even_{\Lambda_\even, \xi_\even}(A)$.
To prove (\ref{equivfin}), first we 
start with the finite volume $\nabla\phi$-Gibbs measure $\mu_{\Lambda,\xi}$. Then we construct a  finite volume $\phi$-Gibbs measure $\nu_{\Lambda,\psi}$ using the map $K_{\Lambda}^\psi$ defined in Remark \ref{equivzd1}. Next we restrict $\nu_{\Lambda,\psi}$ to the even vertices by means of Lemma \ref{nuodd}, and then we pass to the finite volume $\nabla\phi$-Gibbs measure $\mu^\even_{\Lambda_\even,\xi_\even}$ by applying the map $J^{\even,\xi}_{\Lambda_\even}$ defined in Remark \ref{equived}.

The details in the derivation of (\ref{equivfin}) follow below.

Let $\xi\in\chi$. Fixing $\psi(0)\in\RR$, for all $A\in {{\cal F}}_\bzd$ we have by Remark~\ref{equivzd1} that
\begin{equation}
\label{equiv3}
\mu_{\Lambda,\xi}(A)=\ey_{\nu_{\Lambda,\psi}}({1_A\circ {K}^\psi_{\Lambda}}),~\mbox{with}~\psi~\mbox{given as in}~(\ref{19})~\mbox{by}~\psi(x):=\sum_{b\in\C_{0,x}}\xi(b)+\psi(0),\,x\in\zd.
\end{equation}
For all $B\in {\cal F}_{\ed}$ and $\Lambda_\even$ finite sets in $\ed$ with $\ed\setminus\Lambda_\even$ connected, we have by Remark~\ref{equived}
\begin{equation}
\label{equiv3'}
\nu^\even_{\Lambda_\even,{\psi}_\even}(B)=\ey_{\mu^\even_{\Lambda_\even,{\xi}_\even }}(1_B\circ {{J}}^{\even,\xi}_{\Lambda_\even}) ,~\mbox{with}~\xi_\even(b):=\nabla\psi(b),~b\in\bed.
\end{equation}
Let $A\in {{\cal F}}_{\bed}\subset {\cal F}_{\bzd}$; then by using Lemma~\ref{nuodd}, (\ref{equiv3}) and (\ref{equiv3'}), we have for every $\xi\in\chi$ such that  $\xi_\even\in\chi_\even$ (recall Remark~\ref{chitochie})
\begin{multline}
\label{mutomue}
\mu_{\Lambda,\xi}(A)=\ey_{\nu_{\Lambda,\psi}}(1_A\circ {K}^\psi_{\Lambda})=\nu^\even_{\Lambda_\even,\psi_\even}((K^\psi_{\Lambda})^{-1}(A))
=\ey_{\mu^\even_{\Lambda_\even,\xi_\even}}(1_{(K^\psi_{\Lambda})^{-1}(A)}\circ {{J}}^{\even,\xi}_{\Lambda_\even})\\
=\mu^\even_{\Lambda_\even,\xi_\even}(A) ,
\end{multline}
where for the last equality we used the fact that $1_{(K^\psi_{\Lambda})^{-1}(A)}\circ {{J}}^{\even,\xi}_{\Lambda_\even}=1_A$.
\endpf
The following statement is a consequence of the Markov property of the Gibbs measures.
\begin{lemma}
\label{led}{\bf(Conditional of $\nabla\phi$-Gibbs measure on $\bed$)}
Let $G$ be a ${\cal F}_{\bzd}$-measurable and bounded function. Then for all $\mu\in {\cal G}_\beta(H)$ and all $\xi\in\chi$, we have
\begin{equation}
\label{eqninf}
\ey_{\mu}\left(G|{\cal F}_{\bed}\right)(\xi)=\int_{\RR^{\zd}} G(\nabla\phi)\prod_{x\in\od}\nu_{x,\psi}(\rmd\phi(x))\delta_{\psi}(d\phi_{\ed}),
\end{equation}
where we use $\nu_{x, \psi}$ to denote $\nu_{\Lambda, \psi}$ with $\Lambda=\{x\}$
 and $\psi$ is given by $\psi(x):=\sum_{b\in\C^\even_{0,x}}\xi_\even(b)+\psi(0),\,x\in\zd$, for a fixed $\psi(0)\in\RR$ and with $\xi_\even$ given as in Remark \ref{chitochie}. 
\end{lemma}
\pf
It is enough to prove (\ref{eqninf}) for bounded functions $G$ depending on finitely many coordinates. Fix such a $G$ arbitrarily. Note first that the right-hand side of (\ref{eqninf}) is ${\cal F}_{\bed}$-measurable and depends only on the even gradients, as proved in Corollary \ref{tildemu} below. Therefore, to show (\ref{eqninf}) we only need to prove that for any ${\cal F}_{\bed}$-measurable and bounded function $F$ depending on finitely many coordinates in $\bed$, we have
$$\int_{\chi}F({\nabla\phi}_\even) G(\nabla\phi)\mu(\rmd\nabla\phi)=\int_\chi F({\nabla\psi}_\even)\bigg[\int_{\RR^{\zd}} G(\nabla\phi)\prod_{x\in\od}\nu_{x,\psi}(\rmd\phi(x))\delta_{\psi}(d\phi_{\ed})\bigg]\mu(\rmd\nabla\psi).$$ 
Take now an arbitrarily fixed ${\cal F}_{\bed}$-measurable and bounded $F$, depending on finitely many coordinates in $\bed$. For $n\in\N$ let ${\cal S}^d_{n}=\{x\in\zd:||x||\le n\}$ such that $F$ is ${\cal F}_{(\overline{{{\cal S}}^d_{n}})^*}$-measurable and let $\Lambda_n:={\cal S}^d_n\cap\od$. Then from (\ref{dlrgrad}) we have
\begin{equation}
\label{dlr1}
\int_{\chi}F({\eta}_\even) G(\eta)\mu(\rmd\eta)=\int_\chi\mu(d\xi)\int_{\chi_{ \left({{\bar{\Lambda}}_n}\right)^* ,\xi}} \mu_{\Lambda_n,\xi}(d\eta)F(\eta_\even)G(\eta).
\end{equation}
Using Remark \ref{equivzd}, we switch now from the finite $\nabla\phi$-Gibbs measure $\mu_{\Lambda_n,\xi}$ to the corresponding finite $\phi$-Gibbs measure $\nu_{\Lambda_n,\psi}$. Then
\begin{eqnarray}
\label{dlr2}
\int_{\chi_{ \left({{\bar{\Lambda}}_n}\right)^* ,\xi}}\mu_{\Lambda_n,\xi}(d\eta)F(\eta_\even)G(\eta)&=&\int_{\RR^{\zd}}F({\nabla\phi}_\even) G(\nabla\phi)\prod_{x\in\Lambda_n}\nu_{x,\psi}(\rmd\phi(x))\delta_{\psi}(d\phi_{\zd\setminus\Lambda_n})\nonumber\\
&=& F({\nabla\psi}_\even)\int_{\RR^{\zd}}G(\nabla\phi)\prod_{x\in\Lambda_n}\nu_{x,\psi}(\rmd\phi(x))\delta_{\psi}(d\phi_{\zd\setminus\Lambda_n}),
\end{eqnarray}
as $F$ only depends on the even gradients. Since by the Kolmogorov extension theorem we have
$$\int_{\RR^{\zd}} G(\nabla\phi)\prod_{x\in\od}\nu_{x,\psi}(\rmd\phi(x))\delta_{\psi}(d\phi_{\ed})=\lim_{n\rightarrow\infty}\int_{\RR^{\zd}}G(\nabla\phi)\prod_{x\in\Lambda_n}\nu_{x,\psi}(\rmd\phi(x))\delta_{\psi}(d\phi_{\zd\setminus\Lambda_n}),$$
the statement of the theorem follows now from (\ref{dlr1}), (\ref{dlr2}) and Lebesgue's dominated convergence theorem. 
\endpf
In the next Corollary, we reformulate Lemma \ref{led} to remove the dependence on the height
field $\psi$, and to make it more explicit that everything in the formula for $\ey_{\mu}\left(G|{\cal F}_{\bed}\right)(\xi)$ depends only on the
even gradients.
\begin{cor}
\label{tildemuaa}
Let $k\in I$ be an arbitrarily fixed element in $I$ and let $G$ be a ${\cal F}_{\bzd}$-measurable and bounded function. Then for all $\mu\in {\cal G}_\beta(H)$ and all $\xi\in\chi$, we have with the notations from Remark \ref{constrevham} and from Remark \ref{chitochie}
\begin{equation}
\label{hat1mu}
\ey_{\mu}\left(G|{\cal F}_{\bed}\right)(\xi)=\int G\left(\left(\xi_\even(b)-\phi(x)\right)_{b\in {\mathcal B}(x,k),x\in\od}\right)\prod_{x\in\od}\mu_{x,\xi_\even}^k(\rmd\phi(x)),
\end{equation}
where
\begin{equation}
\label{tildemu}
\mu_{x,\xi_\even}^k(\rmd\phi(x))=\frac{1}{Z_{x,\xi_\even}^k}\exp\left(-\beta\sum_{b\in {\mathcal B}(x,k)}U(\xi_\even(b)-\phi(x))\right)\rmd\phi(x),
\end{equation}
and $Z_{x,\xi_\even}^k $ is the normalizing constant.
\end{cor}
\pf
Note first that for all $i\in I$ and all $x\in\od$, $\nabla_i\phi(x)=\phi(x+e_i)-\phi(x+e_k)-\phi(x)+\phi(x+e_k) = \xi_\even(b)-\phi(x)+\phi(x+e_k)$, with $b\in {\mathcal B}(x,k)$. The statement of the corollary follows now immediately, by making in (\ref{eqninf}) the change of variables $\phi(x)\rightarrow\phi(x)+\phi(x+e_k)$ for all $x\in\od$. 
\endpf

\section{Random Walk Representation Condition}

In this section, we prove that under suitable conditions on the perturbation $g$, the new Hamiltonian $H^\even=(H^\even_{\Lambda_\even, \psi_\even})_{\Lambda_\even\subset\ed,\psi_\even\in\ed}$ induced on $\ed$ and defined in (\ref{eqnW}), is strictly convex. More precisely, we will prove that $H^{\even}$ satisfies the so-called \textit{random walk representation condition} (see Definition \ref{rwd} below). This will allow us to adapt results known for strictly convex potentials, such as uniqueness of ergodic component and decay of covariance, to our non-convex setting. 

Subsection 3.1 contains the main result of this section, Theorem \ref{rw}, in which we prove that under assumption (\ref{tag5}) on $g$, the Hamiltonian $H^{\even}$ satisfies the random walk representation condition. Note that, in contrast to the condition in our previous paper \cite{CDM}, $||g''||_{L^{\infty}(\RR)}$ can be arbitrarily large as long as $||g''||_{L^q(\RR)}$ is small. In subsection 3.2, we present some examples of non-convex potentials which fulfill assumption (\ref{tag5}); our first example is the particular class of potentials treated both in \cite{BK} and in \cite{BS}.

\subsection{Definition and Main Result}
For $i\in I$, \llet
$$D^iF_x(y_1,\ldots,y_{d},y_{-1},\ldots,y_{-d}):=\frac{\partial}{\partial y_i} F_x(y_1,\ldots y_{d},y_{-1},\ldots,y_{-d}).$$
We will next formulate a condition on the
multi-body potential, which we call \textbf{the random walk representation condition}, such
that $F_x$ satisfies this condition, and we will adapt earlier results known
for strictly convex two-body potentials to this setting.
\begin{defn}
\label{rwd}
We say that $H^{\even}$ satisfies \textbf{the random walk representation condition} if there exist $\underline{c},\bar{c}>0$ such that for all $x\in\od$, for all $(\phi(x+e_k))_{x\in\od,k\in I}\in\RR^{\ed}$ and all $i,j\in I$
$$\begin{array}{c}
D^{i,i}F_x((\phi(x+e_k))_{k\in I})=-\sum_{j\in I, j\neq i}D^{i,j}F_x((\phi(x+e_k))_{k\in I})\\[1em]
\underline{c}\le -D^{i,j}F_x((\phi(x+e_k))_{k\in I})\le \bar{c}~\mbox{for}~i\neq j.
\end{array}$$
\end{defn}
\begin{rem}\normalfont
Note that for each $x\in\od$, if $H^\even$ satisfies Definition \ref{rwd}, then $F_x$ is uniformly convex (with respect to the even heights). More precisely, for all $\alpha=(\alpha_1,\ldots \alpha_{2d})\in\RR^{2d}$ we have
$$\underline{c}\sum_{i,j\in I,i\neq j}(\alpha_i-\alpha_j)^2\le \sum_{i,j\in I}\alpha_i\alpha_jD^{i,j}F_x((\phi(x+e_k))_{k\in I})\le \bar{c}\sum_{i,j\in I,i\neq j}(\alpha_i-\alpha_j)^2.$$
\end{rem}
\begin{rem}\normalfont
\label{rwrepres}
Potentials satisfying the \textit{random walk representation condition} fulfill the \textit{random walk representation} which is explained, 
for example, in \cite{DGI} or \cite{FSL}. For two-body gradient interactions which are uniformly convex with respect to heights, the random walk representation gives
an extremely useful representation of the covariance matrix, with respect to the measure $\mu_{\Lambda, \xi}$, in terms of the
Green function of a specific random walk.
\end{rem}

The main result of this section is:
\begin{thm}[Random Walk Representation Condition]
\label{rw}
Let $U\in C^2(\RR)$ be such that it satisfies (\ref{tag2}). We also assume that $V,g\in C^2(\RR)$ satisfy (\ref{vc}). Then, if for some $q\ge 1, g''$ satisfies (\ref{tag5}), more precisely, if
$$\beta^{\frac{1}{2q}}||g''||_{L^q(\RR)}<\frac{(C_1)^{\frac{3}{2}}}{2C_2^{\frac{q+1}{2q}}\left(2d\right)^{\frac{1}{2q}}},$$
then there exist $\underline{c},\bar{c}>0$ such that $H^{\even}$ satisfies the random walk representation condition. 
\end{thm}
\begin{rem}
\normalfont
The main idea behind the proof of Theorem \ref{rw} is that one can gain convexity by one-step integration,
which is possible if $||g''||_{L^q(\RR)}$ is sufficiently small compared to $V''$.

What is crucial as regards the bounds $\underline{c},\bar{c}$, is that they are uniform in $x\in\od$ and that they are independent of the possible values of $\phi_\even\in\ed$. This is necessary for us to adapt the arguments known for uniformly strictly convex potentials with two-body interaction to our setting of a generalized random walk representation condition for multi-body potentials.

Note that we only need $||g''||_{L^q(\RR)}$ to be small for the lower bound $\underline{c}$, as the upper bound $\bar{c}$ only 
requires the perturbation to be finite, not small.
\end{rem}
The first step in proving Theorem \ref{rw} is to prove the following lemma
\begin{lemma} 
\label{cov1}
Suppose $x\in\od$. Then for all $j\in I$, we have
\begin{equation}
\begin{array}{rcl}
\label{fee0}
D^{j}F_{x}((\phi(x+e_k))_{k\in I}))&=&-\sum_{i\in I,i\neq j}D^{i}F_{x}((\phi(x+e_k))_{k\in I})),\\[1em]
D^{j,j}F_{x}((\phi(x+e_k))_{k\in I})&=&-\sum_{i\in I,i\neq j}D^{i,j}F_{x}((\phi(x+e_k))_{k\in I}),
\end{array}
\end{equation}
and for all $i\in I, i\neq j$ 
\begin{equation}
\label{feff}
D^{i,j}F_{x}((\phi(x+e_k))_{k\in I})=-4\beta^2\cov_{\nu_{x,\psi_\phi}}\left(U'(\nabla_{i}\phi(x)),U'(\nabla_{j}\phi(x))\right),
\end{equation}
where
$\nu_{x,\psi_\phi}$ is as defined in Lemma \ref{led}, with boundary condition $\psi_\phi(y):=\phi(y)$ for $y\ne x$, and $\ey_{\nu_{x,\psi_\phi}}$ and $\cov_{\nu_{x,\psi_\phi}}$ are respectively the expectation and the covariance with respect to the measure $\nu_{x,\psi_\phi}$.
\end{lemma}
\pf
Let $a=(a_1,a_2,\ldots a_{2d})\in\RR^{2d}.$ Since $F_x(a_1,... a_{2d}) =
F_x(a_1+t,..., a_{2d}+t)~\mbox{for all}~t>0,$
differentiating with respect to $t$ in it, gives the first identity in (\ref{fee0}). The second assertion in (\ref{fee0}) follows from the first, by differentiation. By differentiating now with respect to $\phi(x+e_i)$ and $\phi(x+e_j)$ in the formula for $F_x$, we have for all $i,j\in I, i\neq j$
\begin{equation}
\label{wij}
D^{i,j}F_x((\phi(x+e_k))_{k\in I})=-4\beta^2\cov_{\nu_{x,\psi_\phi}}\left( U'(\nabla_{i}\phi(x)),U'(\nabla_{j}\phi(x))\right).
\end{equation}
\endpf
The next lemma follows by Taylor expansion and will be needed for the proof of Theorem~\ref{rw}:
\begin{lemma} [Representation of Covariances]
\label{covest}

For all $L^2$-functions $F,G\in C^1(\RR;\RR)$ with bounded derivatives and for all measures $\nu\in P(\RR)$, we have
\begin{eqnarray*}
\cov_{\nu} (F,G)&=&\frac{1}{2}\iint\left[F(\phi)-F(\psi)\right]\left[G(\phi)-G(\psi)\right]\nu(\rmd\phi)\nu(\rmd\psi)\nonumber\\
&=&\frac{1}{2}\iint\left[(\phi-\psi)I F(\phi,\psi)\right]\left[(\phi-\psi) I G(\phi,\psi)\right]\nu(\rmd\phi)\nu(\rmd\psi),\nonumber
\end{eqnarray*}
where we denote by
$$I F(\phi,\psi):=\int_0^1 F'\left(\psi+t(\phi-\psi)\right)\rmd t,~~I G(\phi,\psi):=\int_0^1 G'\left(\psi+s(\phi-\psi)\right)\rmd s.$$
\end{lemma}
\begin {rem} [Scaling Argument]
\label{scal}
\normalfont
A simple scaling argument shows that it suffices to prove Theorem \ref{rw} for 
\begin{eqnarray}
\label{11}
\beta=1, C_1=1.
\end{eqnarray}
Indeed, suppose that the result is true for $\beta=1$ and $C_1=1$. Given $\beta$, $V$ and $g$ which satisfy (\ref{vc}) and (\ref{tag5}), we define
$$\tilde{U}(s)=\tilde{V}(s)+\tilde{g}(s),~\mbox{where}~\tilde{V}(s)=\beta V\left(\frac{s}{\sqrt{\beta C_1}}\right),~\tilde{g}(s)=\beta g\left(\frac{s}{\sqrt{\beta C_1}}\right).$$
Then
$$1\le (\tilde{V})''\le\frac{C_2}{C_1},~-\frac{C_0}{C_1}\le (\tilde{g})''\le 0,||\tilde{g})''||_{L^q(\RR)}={(\beta C_1)}^{\frac{1}{2q}}\frac{1}{C_1}||g''||_{L^q(\RR)},||\tilde{g})'||_{L^2(\RR)}={(\beta^3/C_1)}^{\frac{1}{4}}||g'||_{L^2(\RR)}.$$
Hence $\tilde{V}$, $\tilde{g}$ satisfy the assumptions of Theorem~\ref{rw} with $\beta=1$ and $C_1=1$. 
On the other hand, the change of variables $\tilde{\phi}(x)=\sqrt{\beta C_1}\phi(x)$ yields $\tilde{U}\left(\nabla_i\tilde{\phi}(x)\right)=\beta U(\nabla_i\phi(x))$ and thus
\begin{multline*}
\tilde{F}_x((\tilde{\phi}(x+e_i))_{i\in I}):= -\log\int_{\RR} e^{-2\sum_{i\in I} \tilde{U}(\btd_{i}\tilde{\phi}(x))}\rmd\tilde{\phi}(x)\\
=-\frac{\log{\beta C_1}}{2}-\log\int_{\RR} e^{-2\beta\sum_{i\in I} U(\btd_{i}\phi(x))}\rmd\phi(x)=-\frac{\log{\beta C_1}}{2}+F_x((\phi(x+e_i))_{i\in I}).
\end{multline*}
\end{rem}
~\\
\noindent\textbf{Proof of Theorem~\ref{rw}}
From Definition~\ref{rwd} and Lemma~\ref{cov1} it follows that, in order to prove that the random walk representation condition holds for $H^{\even}$, all we need is to show that there exist $c_l,c_u>0$ such that for all $i,j\in I, i\neq j$, and uniformly in $x$ and $\psi$ 
\begin{equation}
\label{lowupcov}
c_l\le\cov_{\nu_{x,\psi}}\left(U'(\nabla_{i}\phi(x)),U'(\nabla_{j}\phi(x))\right)\le c_u.
\end{equation}
Recall that we have $U = V+ g$, where $1\le V'' \le C_2$ and therefore we can split the initial covariance term into four resulting covariance terms. More precisely, we have
$$cov_{\nu_{x,\psi}}(U'_i,U'_j)=cov_{\nu_{x,\psi}}(V'_i,V'_j)+cov_{\nu_{x,\psi}}(V'_i,g'_j)+cov_{\nu_{x,\psi}}(V'_j,g'_i)+cov_{\nu_{x,\psi}}(g'_i,g'_j),$$
where for convenience of notation we denote by $$cov_{\nu_{x,\psi}}(U'_i,U'_j):=cov_{\nu_{x,\psi}}\left(U'(\nabla_{i}\phi(x)),U'(\nabla_{j}\phi(x))\right),\ldots,cov_{\nu_{x,\psi}}(g'_i,g'_j):=cov_{\nu_{x,\psi}}\left(g'(\nabla_{i}\phi(x)),g'(\nabla_{j}\phi(x))\right).$$
We will first show in (\ref{vi}), (\ref{gi}) and (\ref{vx3a}) below, by means of Lemma \ref{covest}, that the resulting $\cov_{\nu_{x,\psi}}(V'_i,V'_j)$ and $\cov_{\nu_{x,\psi}}(g'_i,g'_j)$ terms are positive and that the resulting $\cov_{\nu_{x,\psi}}(g'_i,V'_j)$ and  $\cov_{\nu_{x,\psi}}(g'_j,V'_i)$ terms are negative. We will then obtain lower and upper bound estimates for the $\cov_{\nu_{x,\psi}}(V'_i,V'_j)$ terms, and upper bound estimates for the  $\cov_{\nu_{x,\psi}}(g'_i,g'_j)$ and the $-\cov_{\nu_{x,\psi}}(g'_i,V'_j)$ and $-\cov_{\nu_{x,\psi}}(g'_j,V'_i)$ terms. These bounds will determine the conditions on the perturbation $g''$ such that (\ref{lowupcov}) holds. To estimate an arbitrary $\cov_{\nu_{x,\psi}}(V'_i,V'_j)$ term, we will bound it in (\ref{vi}) from above and below by bounds proportional to  $\cov_{\nu_{x,\psi}}(\phi,V'_j)$. To estimate an arbitrary $\cov_{\nu_{x,\psi}}(g'_i,g'_j)$ term, we will bound the respective term in (\ref{gi}) from above by a bound proportional to  $\cov_{\nu_{x,\psi}}(\phi,V'_j)$. To estimate an arbitrary $-\cov_{\nu_{x,\psi}}(V'_j,g'_i)$ term, we will first express it in (\ref{altvar}) in terms of $\cov_{\nu_{x,\psi}}(\phi,V'_j)$ and $\var_{\nu_{x,\psi}}(g'_i)$; the $\var_{\nu_{x,\psi}}(g'_i)$ term will then also be bound in (\ref{vx3}) from above by a bound proportional to $\cov_{\nu_{x,\psi}}(\phi,V'_j)$. We will then proceed to find upper and lower bounds for the $\cov_{\nu_{x,\psi}}(\phi,V'_j)$ terms. The upper bound will be derived in (\ref{newvga}) by means of (\ref{vx1}), (\ref{vx2}) and (\ref{l1lbound}), and the lower bound will be derived in (\ref{lbound}) by means of (\ref{vx2}). The explicit computations follow. 

Fix $x\in\zd$ and $i,j\in I,i\neq j,$ arbitrarily. We will next check which covariance terms are positive and which are negative. Using Lemma~\ref{covest} for $V'(\nabla_i\phi(x))$ and $V'(\nabla_j\phi(x))$, we see that 
\begin{multline*}
\cov_{\nu_{x,\psi}} (V'(\nabla_{i}\phi(x)), V' (\nabla_{j}\phi(x)))=
\frac{1}{2}\iint(\phi(x)-\psi(x))^2\int_0^1 V''\left((1-t)\psi(x)-\phi(x+e_i)+t\phi(x)\right)dt\\
\int_0^1 V''\left((1-s)\psi(x)-\phi(x+e_j)+s\phi(x)\right)\rmd s\nu_{x}(\rmd\phi)\nu_{x}(\rmd\psi).
\end{multline*}
By comparing the above equality with the similar one for $\cov_{\nu_{x,\psi}}(\phi(x),V' (\nabla_{j}\phi(x)))$ and with the bound $1 \le V'' \le C_2$, we have for all $i,j\in I$
\begin{equation}
\begin{array}{l}
\label{vi}
\cov_{\nu_{x,\psi}} (V'(\nabla_{i}\phi(x)), V' (\nabla_{j}\phi(x)))\ge \cov_{\nu_{x,\psi}}(\phi(x),V' (\nabla_{j}\phi(x)))\ge \var_{\nu_{x,\psi}}(\phi(x))\ge 0,\\[1em]
\cov_{\nu_{x,\psi}} (V'(\nabla_{i}\phi(x)), V' (\nabla_{j}\phi(x)))\le C_2\cov_{\nu_{x,\psi}}(\phi(x),V' (\nabla_{j}\phi(x))).
\end{array}
\end{equation}
Since $-C_0\le g''\le 0$, by similar reasoning
\begin{equation}
\label{gi}
0\le\cov_{\nu_{x,\psi}} (g'(\nabla_{i}\phi(x)), g'(\nabla_{j}\phi(x)))\le C_0^2\var_{\nu_{x,\psi}}(\phi(x))\le C_0^2\cov_{\nu_{x,\psi}}(\phi(x),V' (\nabla_{j}\phi(x))),
\end{equation} 
and
\begin{equation}
\label{vx3a}
-C_0\cov_{\nu_{x,\psi}}(\phi(x),V'(\nabla_{j}\phi(x))) \le\cov_{\nu_{x,\psi}} (V'(\nabla_{j}\phi(x)), g'(\nabla_{i}\phi(x)))<0.
\end{equation}
Given (\ref{vi}), (\ref{gi}) and (\ref{vx3a}), we have the following upper and lower bounds for $\cov_{\nu_{x,\psi}}(U',U')$
\begin{multline}
\label{covvg}
 \cov_{\nu_{x,\psi}} (\phi(x), V' (\nabla_{j}\phi(x)))+\cov_{\nu_{x,\psi}} (g'(\nabla_{j}\phi(x)), V' (\nabla_{i}\phi(x)))+\cov_{\nu_{x,\psi}} (g'(\nabla_{i}\phi(x)), V' (\nabla_{j}\phi(x)))\\ \le\cov_{\nu_{x,\psi}}\left(U'(\nabla_{i}\phi(x)),U'(\nabla_{j}\phi(x))\right)\le\left(C_2+C_0^2\right)\cov_{\nu_{x,\psi}} (\phi(x), V' (\nabla_{j}\phi(x))).
\end{multline}
Of more importance are the lower bound estimates, as they will determine the conditions on our perturbation $g''$ which give us convexity after the one-step integration. We will next get a lower bound for the $\cov_{\nu_{x,\psi}} (g'_i, V'_j)$ terms in (\ref{covvg}), which shows that the upper and lower bounds in (\ref{covvg}) are all in terms of $\cov_{\nu_{x,\psi}} (\phi, V'_j)$. Using (\ref{vx3a}), the Cauchy-Schwarz inequality and (\ref{vi}), we have
\begin{multline}
\label{altvar}
0\le-\cov_{\nu_{x,\psi}} (V'(\nabla_{j}\phi(x)), g'(\nabla_{i}\phi(x)))\le\sqrt{\var_{\nu_{x,\psi}}(V'(\nabla_{j}\phi(x)))}\sqrt{\var_{\nu_{x,\psi}}(g'(\nabla_{i}\phi(x)))}\\
\le \sqrt{C_2\cov_{\nu_{x,\psi}}(\phi(x),V'(\nabla_{j}\phi(x)))}\sqrt{\var_{\nu_{x,\psi}}(g'(\nabla_{i}\phi(x)))}.
\end{multline}
Let now $q\ge 1$ be arbitrarily fixed. By Lemma~\ref{covest} and Jensen's inequality, we get
\begin{eqnarray}
\label{vx3}
\lefteqn{\var_{\nu_{x,\psi}}(g'(\nabla_{i}\phi(x)))}\nonumber\\
&=&\frac{1}{2}\iint(\phi(x)-\psi(x))^2\left[\int_0^1 g''\left(\psi(x)-\phi(x+e_i)+t(\phi(x)-\psi(x))\right)\rmd t\right]^2\nu_{x}(\rmd\phi)\nu_{x}(\rmd\psi)\nonumber\\
&\le&\frac{1}{2}\iint(\phi(x)-\psi(x))^2\left[\int_0^1 \left|g''\left(\psi(x)-\phi(x+e_i)+t(\phi(x)-\psi(x))\right)\right|^q\rmd t\right]^{\frac{2}{q}}\nu_{x}(\rmd\phi)\nu_{x}(\rmd\psi)\nonumber\\
&=&\frac{1}{2}\iint|\phi(x)-\psi(x)|^{2-2/q}\left[\int_{\psi(x)-\phi(x+e_i)}^{\phi(x)-\phi(x+e_i)} \left|g''\left(s\right)\right|^q\rmd s\right]^{\frac{2}{q}}~\nu_{x}(\rmd\phi)\nu_{x}(\rmd\psi)\nonumber\\
&\le&\frac{1}{2}||g''||^2_{L^q(\RR)}\iint|\phi(x)-\psi(x)|^{2-2/q}~\nu_{x}(\rmd\phi)\nu_{x}(\rmd\psi)\le \frac{1}{2^{\frac{1}{q}}}||g''||^2_{L^q(\RR)}\left[\var_{\nu_{x,\psi}}(\phi(x))\right]^{\frac{q-1}{q}}\nonumber\\
&\le& \frac{1}{2^{\frac{1}{q}}}||g''||^2_{L^q(\RR)}\left[\cov_{\nu_{x,\psi}}(\phi(x),V'(\nabla_{j}\phi(x)))\right]^{\frac{q-1}{q}}.
\end{eqnarray}
where for the second equality we made the change of variable $s=\psi(x)-\phi(x+e_i)+t(\phi(x)-\psi(x))$, in the penultimate inequality we used Lemma~\ref{covest} and for the last inequality we used (\ref{vi}). The lower bound in (\ref{covvg}) becomes by (\ref{vx3})
\begin{multline}
\label{lbound1}
\lefteqn{\cov_{\nu_{x,\psi}}\left(U'(\nabla_{i}\phi(x)),U'(\nabla_{j}\phi(x))\right)}\\
\ge\left[\cov_{\nu_{x,\psi}}(\phi(x),V'(\nabla_{j}\phi(x)))\right]^{\frac{2q-1}{2q}}\left[\left[\cov_{\nu_{x,\psi}}(\phi(x),V'(\nabla_{j}\phi(x)))\right]^{\frac{1}{2q}}-2^{(2q-1)/{2q}}\sqrt{C_2}||g''||_{L^q(\RR)}\right].
\end{multline}
We now proceed to find upper and lower bounds for $\cov_{\nu_{x,\psi}}(\phi(x),V'(\nabla_{j}\phi(x)))$. From (\ref{vi}), we have by repeated application
\begin{equation}
\label{vx1}
\cov_{\nu_{x,\psi}}(\phi(x),V' (\nabla_{j}\phi(x)))\le\frac{1}{2d}\cov_{\nu_{x,\psi}}\left(V'(\nabla_j\phi(x)),\sum_{i\in I} V'(\nabla_{i}\phi(x)))\right).
\end{equation}
Recall now that
\begin{multline*}
\cov_{\nu_{x,\psi}}\left(V'(\nabla_j\phi(x)),\sum_{i\in I} V'(\nabla_{i}\phi(x))\right) =\frac{1}{Z_{x,\psi}}\int V'(\nabla_j\phi(x))\left(\sum_{i\in I} V'(\nabla_{i}\phi(x))\right)e^{-2 H_x(\phi)}d\phi(x)\\
-\left[\frac{1}{Z_{x,\psi}}\int V'(\nabla_j\phi(x))e^{-2 H_x(\phi)}d\phi(x)\right]\left[\frac{1}{Z_{x,\psi}}\int \left(\sum_{i\in I} V'(\nabla_{i}\phi(x))\right)e^{-2H_x(\phi)}d\phi(x)\right],
\end{multline*}
where $Z_{x,\psi}$ is the normalizing constant and $H_x(\phi)$ has been defined in (\ref{hxeqn}).
Using integration by parts in the above, we have
\begin{multline}
\label{vx2}
\cov_{\nu_{x,\psi}}\left(V'(\nabla_j\phi(x)),\sum_{i\in I} V'(\nabla_{i}\phi(x))\right)=\frac{1}{2}\ey_{\nu_{x,\psi}}\left(V''(\nabla_{j}\phi(x))\right)\\-\cov_{\nu_{x,\psi}}\left(V'(\nabla_{j}\phi(x)),\sum_{i\in I}g'(\nabla_{i}\phi(x))\right)\le \frac{C_2}{2}-\cov_{\nu_{x,\psi}}\left(V'(\nabla_{j}\phi(x)),\sum_{i\in I}g'(\nabla_{i}\phi(x))\right).
\end{multline}
 From (\ref{vx1}), (\ref{vx2}) and (\ref{vx3}), we now get the upper bound
$$\cov_{\nu_{x,\psi}}(\phi(x),V'(\nabla_{j}\phi(x)))\le\frac{C_2}{4d}+ \frac{\sqrt{C_2}}{2^{(2q+1)/2q}d}||g''||_{L^q(\RR)}\left[\cov_{\nu_{x,\psi}}(\phi(x),V'(\nabla_{j}\phi(x)))\right]^{\frac{2q-1}{2q}},$$
which is equivalent to 
\begin{equation}
\label{l1lbound}
\left[\cov_{\nu_{x,\psi}}(\phi(x),V'(\nabla_{j}\phi(x)))\right]^{\frac{2q-1}{2q}}\left[ \left[\cov_{\nu_{x,\psi}}(\phi(x),V'(\nabla_{j}\phi(x)))\right]^{\frac{1}{2q}}-b\right]\le a,
\end{equation}
where $a=\frac{C_2}{4d}$ and $b=\frac{\sqrt{C_2}}{2^{(2q+1)/2q}d}||g''||_{L^q(\RR)}$. Depending on if $\left[\cov_{\nu_{x,\psi}}(\phi(x),V'(\nabla_{j}\phi(x)))\right]^{\frac{1}{2q}}\le b$ or $\ge b$, (\ref{vi}) combined with simple arithmetic in the above inequality gives 
\begin{multline}
\label{newvga}
\tau^2_{x,\psi}:=\var_{\nu_{x,\psi}}(\phi(x)))\le\cov_{\nu_{x,\psi}}(\phi(x),V'(\nabla_{j}\phi(x)))\le\max\left[b^{2q},\left(\frac{a}{b^{\frac{2q-1}{2q}}}+b\right)^{2q} \right]=\left(\frac{a}{b^{\frac{2q-1}{2q}}}+b\right)^{2q}.
\end{multline}
The upper bound on $\cov_{\nu_{x,\psi}}\left(U'(\nabla_{i}\phi(x)),U'(\nabla_{j}\phi(x))\right)$ follows now from (\ref{covvg}) and (\ref{newvga}).
To find a lower bound, note now that from (\ref{vi}) we get
$$\cov_{\nu_{x,\psi}}(\phi(x),V' (\nabla_{j}\phi(x)))\ge\frac{1}{2dC_2}\cov_{\nu_{x,\psi}}\left(V'(\nabla_j\phi(x)),\sum_{i\in I} V'(\nabla_{i}\phi(x))\right).$$
By using (\ref{vx2}) and (\ref{vx3a}), we have
\begin{equation}
\label{lbound}
\cov_{\nu_{x,\psi}}(\phi(x),V' (\nabla_{j}\phi(x)))\ge\frac{1}{4d C_2}.
\end{equation}
From (\ref{lbound}) and (\ref{lbound1}), the lower bound becomes
$$\cov_{\nu_{x,\psi}}\left(U'(\nabla_{i}\phi(x)),U'(\nabla_{j}\phi(x))\right)\ge\frac{1}{\left(4d C_2\right)^{\frac{2q-1}{2q}}}\left[\frac{1}{\left(4d C_2\right]^{\frac{1}{2q}}}-\frac{2\sqrt{C_2}||g''||_{L^q(\RR)}}{2^{\frac{1}{2q}}}\right].$$
 To summarize, we obtain the following upper and lower bounds, uniform with respect to $x$ and $\psi$
\begin{equation}
\label{loboundcov}
c_l=\frac{1}{\left(4d C_2\beta\right)^{\frac{2q-1}{2q}}} \epsilon\le\cov_{\nu_{x,\psi}}\left(U'(\nabla_{i}\phi(x)),U'(\nabla_{j}\phi(x))\right)\le \left(C_2+C_0^2\right)\left(\frac{a}{b^{\frac{2q-1}{2q}}}+b\right)^{2q}=c_u,
\end{equation}
for $\epsilon=\frac{1}{\left(4d C_2\right)^{\frac{1}{2q}}}-\frac{2\sqrt{C_2}||g''||_{L^q(\RR)}}{2^{\frac{1}{2q}}}>0$ by (\ref{tag5}).  
\endpf
\begin{rem}\normalfont 
Another possible condition, (\ref{g'l2}), is obtained if we use Lemma \ref{l1norm} below to replace (\ref{vx3}) by
$$\var_{\nu_{x,\psi}}(g'(\nabla_{i}\phi(x)))\le\eyb_{\nu_{x,\psi}}\left((g'(\nabla_{i}\phi(x)))^2\right)\le 2\sqrt{\beta dC_2}||g'||^2_{L^2(\RR)}.$$
\end{rem}
\begin{lemma}
\label{l1norm}
If $h\in L^1(\RR)$, then we have
$$\left|\ey_{\nu_{x,\psi}} (h)\right|\le 2\sqrt{d\beta C_2}||h||_{L^1(\RR)}.$$
\end{lemma}
\pf
Using integration by parts and Cauchy-Schwarz, we have
\begin{eqnarray*}
\left|\ey_{\nu_{x,\psi}}(h)\right|&=&\left|\eyb_{\nu_{x,\psi}}\left(\frac{\partial}{\partial y}\left(\int_{-\infty}^y h(z)\rmd z\right)\right)\right|
=2\beta\left|\eyb_{\nu_{x,\psi}}\left(H_x'(y)\left(\int_{-\infty}^y h(z)\rmd z\right)\right)\right|\nonumber\\
&\le& 2\beta\left[\eyb_{\nu_{x,\psi}}\left((H_x')^2\right)\right]^{1/2}\left[\eyb_{\nu_{x,\psi}}\left(\left(\int_{-\infty}^y h(z)\rmd z\right)^2\right)\right]^{1/2}\nonumber\\
&=&\sqrt{2\beta}\left[\eyb_{\nu_{x,\psi}}(H_x'')\right]^{1/2}\left[\eyb_{\nu_{x,\psi}}\left(\int_{-\infty}^y h(z)\rmd z\right)^2\right]^{1/2}\le 2\sqrt{d\beta C_2}||h||_{L^1(\RR)}.
\end{eqnarray*}
Note that we also used property (\ref{vc}) in the above formula.
\endpf

\begin{rem}\normalfont
\label{strictconv}
Note that if we consider the case where $U$ is strictly convex with $C_1\le U''\le C_2$ (that is $U=V$ and $g=0$), in view of (\ref{vi}) and (\ref{l1lbound}), the one step integration preserves the strict convexity of the induced Hamiltonian as
$$\frac{C_1^2}{4d\beta C_2}\le\cov_{\nu_{x,\psi}}\left(U'(\nabla_{i}\phi(x)),U'(\nabla_{j}\phi(x))\right)\le\frac{C_2^2}{4d\beta C_1}.$$
\end{rem}

\begin{rem}\normalfont \textbf{(Perturbation with Compact Support)}
\label{compsup}
Note that we can extend the results from Theorem \ref{rw} to the case where we have a perturbation $g$ such that $g''$ has compact support (see also example (b) below). More precisely, assume that $U=Y+h$, where $U$ 
satisfies (\ref{tag2}), $D_1\le Y''\le D_2$ and $-D_0\le h''\le 0$ on $[a,b]$ and $0<h''<D_3$ on $\RR\setminus [a,b]$, with $a,b\in\RR$ and $h''(a)=h''(b)=0$. Then we just need to replace 
$$C_1:=D_1, C_2:=D_1+D_2,~\mbox{and}~g'':=h''1_{\{h''\le 0\}}.$$ 
A sketch of the argument follows next. Set
$$g(s)=h(s)1_{\{s\in [a,b]\}}+\left[h(b)+h'(b)(s-b)\right]1_{\{s>b\}}+\left[h(a)+h'(a)(s-a)\right]1_{\{s<a\}}$$
and 
$$V(s)=Y(s)+h(s)1_{\{s\notin [a,b]\}}-\left[h(b)+h'(b)(s-b)\right]1_{\{s>b\}}-\left[h_i(a)+h'(a)(s-a)\right]1_{\{s<a\}}.$$
Thus, we have $V,g\in C^2(\RR)$, with $-D_0\le h''(s)=g''(s)\le 0$ for $s\in[a,b]$ and $g''(s)=0$ for $s\in\RR\setminus [a,b]$ and $D_1\le V''(s)=Y''(s)+h''(s)1_{\{s\notin [a,b]\}}\le D_2+D_3$. Note that this procedure can also be extended to the case where $h''$ changes sign more than once. 
\end{rem}

\subsection{Examples}

\begin{enumerate}

\item [(a)] 

 Let $p\in (0,1)$ and $0<k_2<k_1$. Let
$$U(s)=-\log\left(pe^{-k_1\frac{s^2}{2}}+(1-p)e^{-k_2\frac{s^2}{2}}\right).$$
Take $\frac{p}{1-p}>\frac{k_2}{k_1}$ in order that the potential $U$ is non-convex. Let $\beta=1$, $d=2$ and $k_1\gg k_2$. In this particular case, as Christof K\"ulske pointed out to us, we are dealing entirely with sums of Gaussian integrals, so we can compute $\cov_{\nu_{x,\psi}}\left(U'(\nabla_{i}\phi(x)),U'(\nabla_{j}\phi(x))\right)$ directly, which explicit computation is not possible in general; the random walk representation condition holds then if $\frac{p}{1-p}<O\left(\left(\frac{k_2}{k_1}\right)^{1/2}\right)$ (see the Appendix for a sketch of the explicit computations).

This particular example is of independent interest and has been the focus of two other papers in the area (see \cite{BK} and \cite{BS}). For the case $d=2$ and $\beta=1$, it was proved in \cite{BK} that at the critical point $p:=p_c$, such that $\frac{p_c}{1-p_c}=\left(\frac{k_2}{k_1}\right)^{1/4}$, uniqueness of ergodic states is violated for this example of potential $U$ and there are multiple ergodic, invariant $\nabla\phi$-Gibbs measures with zero tilt; the same example is also treated in \cite{BS}, where they prove CLT for the this particular class of potentials in the case of $\nabla\phi$-Gibbs measures with zero tilt.

Note that we can use (\ref{g'l2}) to show that the random walk representation condition holds if $p<O\left(\left(\frac{k_2}{k_1}\right)^{2/3}\right)$. To show this, take $V$ and $g$ even, with $V(0)=0$, $g(0)=0$, and such that
\begin{equation}
\label{3.2a}
V''(s)=\frac{pk_1e^{-k_1\frac{s^2}{2}}+(1-p)k_2e^{-k_2\frac{s^2}{2}}}{pe^{-k_1\frac{s^2}{2}}+(1-p)e^{-k_2\frac{s^2}{2}}}, g''(s)=-\frac{p(1-p)(k_1-k_2)^2s^2}{p^2e^{-(k_1-k_2)\frac{s^2}{2}}+2p(1-p)+(1-p)^2 e^{(k_1-k_2)\frac{s^2}{2}}}.
\end{equation} 
Then
$$k_2\le V''(s)\le p k_1+(1-p)k_2,~||g'(s)||_{L^2(\RR)}\le O\left(\frac{p}{1-p}(k_1-k_2)^{1/4}\right),$$
$$\frac{p}{1-p}(k_1-k_2)^{1/4}\le O\left(\frac{(k_2)^{3/2}}{(pk_1+(1-p)k_2)^{5/4}}\right)=O\left(\frac{(k_2)^{3/2}}{(pk_1)^{5/4}}\right).$$

\begin{figure}[!h]
~\hfill
\begin{minipage}[t]{2in}
\input{ExampleA1.tex.inc}
\caption{Example (a)}
\end{minipage}
\hfill
\begin{minipage}[t]{2in}
\input{ExampleB1.tex.inc}
\caption{Example (b)}
\end{minipage}
\hfill~
\end{figure}

\item [(b)] $U(s)=s^2+a-\log(s^2+a),~~\mbox{where}~~0<a<1$. Let  $0<\beta<\frac{a}{4\sqrt{2}d\left(2+\frac{2}{25a}\right)^2}$. This 
example is interesting, as it has two global minima.

Then, using the notation from Remark \ref{compsup}, take $Y(s)=s^2$ and $h(s)=-\log(s^2+a)$. We have $Y''(s)=2$, so $D_1=D_2=2$; also $h''(s)=2\frac{s^2-a}{(s^2+a)^2}$, with $-\frac{2}{a}\le h''(s)\le 0$ for $s\in [-\sqrt{a},\sqrt{a}]$ and $0<h''(s)\le\frac{2}{25a}$ otherwise. Then $C_0=\frac{2}{a}$, $C_1=2$,$C_2= 2+\frac{2}{25a}$ and $||g''(s)||_{L^1(\RR)}=\frac{2}{\sqrt{a}}$. By using condition (\ref{tag5}) with $q=1$, the random walk representation condition holds.


\end{enumerate}

\section{ Uniqueness of ergodic component}
In this section, we extend the uniqueness of ergodic component result, proved for strictly convex potentials in \cite{FS}, to the class of non-convex potentials $U=V+g$ which satisfy (\ref{tag2}) such $V$ 
and $g$ satisfy (\ref{vc}) and (\ref{tag5}). Note that existence of an ergodic $\mu_u$ is guaranteed for our class of non-convex potentials by Theorem \ref{existerg} below.

The proof of Theorem \ref{ergall} will be done in two steps. First, in subsection 5.1 we will prove the uniqueness of ergodic, shift-invariant $\mu^{\even}_u\in {\cal G}_\even(H^{\even})$ with a given tilt $u\in\RR^d$, when the potentials $F_x$ are of form as defined in (\ref{hameven}) and therefore $H^{\even}$ satisfies the random walk representation condition. For that, we will be adapting earlier
results for two-body potentials under uniformly strictly convex condition, to multi-body potentials satisfying the random walk
representation condition. Then we will use this result combined with Lemma \ref{led} in subsection 5.2, to extend the result to $\mu_u\in {\cal G}_\beta(H)$.

\subsection{Step 1: Uniqueness of ergodic component for $\bed$}
For $x\in\ed$, we define the \textbf{even} shift operators: $\sigma_{x}:\RR^{\ed}\rightarrow\RR^{\ed}$ and $\sigma_{x}:\RR^{\bed}\rightarrow\RR^{\bed}$ similarly as for $x\in\zd$. Then shift-invariance and ergodicity for $\mu^{\even}$ (with respect to $\sigma_x$ for all $x\in\ed$) are defined similarly as for $\mu$. The main result in this section is: 
\begin{thm}
\label{ergeven1}
For every $u\in\RR^d$, there exists at most one $\mu^{\even}_u\in {\cal G}_\even(H^{\even})$, shift-invariant and ergodic with tilt $u$.
\end{thm}
We will prove Theorem \ref{ergeven1} by coupling
techniques. We will follow the same line of argument as in \cite{FS}, by introducing dynamics on the
gradient field which keeps the measure in ${\cal G}_\even(H^\even)$ invariant. Suppose the dynamics of the \textbf{even} height variables $\phi_t=\{\phi_t(y)\}_{y\in\ed}$ are generated by the family of SDEs
\begin{equation}
\label{sde}
\rmd\phi_t(y)=-\sum_{x\in\od,||x-y||=1}\frac{\partial}{\partial\phi(y)}F_x((\phi_t(x+e_i))_{i\in I})\rmd t+\sqrt{2}d W_t(y),~~y\in\ed,
\end{equation}
where for all $x\in\od$, $F_x$ are the functions defined in Lemma \ref{nuodd}, satisfying the properties in Definition \ref{rwd}, and $\{W_t(y),y\in\ed\}$ is a family of independent Brownian motions. Using standard SDE methods and due to the fact that $V''$ is bounded, one can show that equation (\ref{sde}) has a unique solution in $L_r^2$ for some $r>0$. 

We denote by $S_\even$ the class of all shift invariant $\mu\in P_2({\chi}_\even)$ which are 
stationary for the SDE (\ref{sde}) and by ext $S_\even$ those $\mu_\even\in S_\even$ which are ergodic. For each $u\in\RR^d$, we denote by $\left(\ext{{\cal S}_\even}\right)_{u}$ 
the family of all $\mu^\even\in\ext {\cal S}_\even$ such that $\ey_{\mu^\even}(\eta_\even(b))= \langle u, y_b-x_b\rangle$ for all bonds $b=(x_b,y_b)\in\bed$. Note that all translation invariant measures in ${\cal G}_\even(H^\even)$ are stationary under the dynamics (see Proposition 3.1 in \cite{FS}).

The next theorem is a key result in the proof of Theorem \ref{ergeven1}.
\begin{thm}
\label{ergeven}
For every $u\in\RR^d$, there exists at most one $\mu^{\even}_{u}\in \left(\ext{{\cal S}_\even}\right)_{u}$.
\end{thm}
Theorem \ref{ergeven1} now follows from Theorem \ref{ergeven} and Proposition 3.1 in \cite{FS}, which shows that if $\mu^\even_u\in {\cal G}_\even(H^\even)$ is shift-invariant and ergodic, then $\mu^\even_u\in\ext{{\cal S}_\even}$.

The proof of Theorem \ref{ergeven} is based on a coupling lemma, Lemma \ref{2.1} below; a key ingredient for the coupling lemma is a bound
on the distance between two measures evolving under the same dynamics. The main ingredients needed to prove it are Lemma \ref{lsde} below and a non-standard 
ergodic theorem (see (\ref{ergtheo}) below). The deduction of Theorem \ref{ergeven} from the coupling lemma follows the same arguments as the proof of Theorem 2.1 in \cite{FS} and will be omitted.

~

\noindent \textbf{Dynamics}
We will first derive a differential inequality for the difference of two solutions evolving under
the same dynamics, which will be a key ingredient in the proof of the
coupling Lemma \ref{2.1} below.
\begin{lemma}
\label{lsde}
Let $\phi_t$ and $\bar{\phi}_t$ be two solutions for (\ref{sde}), coupled via the same Brownian motion in (\ref{sde}), and set $\tilde{\phi_t}(y):=\phi_t(y)-\bar{\phi_t}(y)$, where $y\in\ed$. Then for every 
finite $\Lambda_\even\subset\ed$, we have
\begin{equation}
\label{bi}
\frac{\partial}{\partial t}\sum_{y\in\Lambda_\even} (\tilde{\phi_t}(y))^2\le -\underline{c}\sum_{b\in\bedl}\left[\nabla\tilde{\phi_t}(b)\right]^2+2\bar{c}\sum_{b\in\partial\bedl}|\phi_t(y_b)| |\nabla\tilde{\phi_t}(b)|.
\end{equation}
\end{lemma}
\pf
The proof of Lemma \ref{lsde} is an adaptation of an earlier result by \cite{FS}, where we replace the uniform strictly convex condition on the two-body potential $V$ with the random walk representation condition on a multi-body potential of gradient type.

Let $y\in\Lambda_\even$. Then from (\ref{sde}), we have
\begin{equation}
\label{newsde}
\frac{\partial}{\partial t}(\tilde{\phi_t}(y))^2=-2\sum_{x\in \odl,||x-y||=1}\bigg[\frac{\partial}{\partial\phi(y)}F_x(({\phi}_t(x+e_i))_{i\in I})-\frac{\partial}{\partial\phi(y)}F_x((\bar{\phi}_t(x+e_i))_{i\in I})\bigg]\tilde{\phi_t}(y).
\end{equation}
By summing now in (\ref{newsde}) over all $y\in\Lambda_\even$ in (\ref{newsde}), we get 
\begin{equation}
\label{3}
\frac{\partial}{\partial t}\sum_{y\in\Lambda_\even} (\tilde{\phi_t}(y))^2=
-2\sum_{x\in\odl}\sum_{\{j\in I|\atop x+e_j\in\Lambda_\even\}}\bigg[D^jF_x((\phi_t(x+e_i))_{i\in I})-D^jF_x((\bar{\phi}_t(x+e_i))_{i\in I})\bigg]\tilde{\phi_t}(x+e_{j}),
\end{equation}
where $\odl=\Lambda\cap\od$ and $\Lambda$ is the associated set to $\Lambda_\even$, as defined in Definition \ref{ascolambda}. To prove (\ref{bi}), we expand now $D^jF_x(({\phi}_t(x+e_i))_{i\in I})$ around $(\bar{\phi}_t(x+e_i))_{i\in I}$ by the Mean Value Theorem to get
\begin{multline}
\label{tay}
D^j F_x(({\phi}_t(x+e_i))_{i\in I})-D^j F_x((\bar{\phi}_t(x+e_i))_{i\in I})\\
=\sum_{k\in I}\tilde{\phi_t}(x+e_k)\int_0^1 D^{j,k}F_x\left((s{\phi}_t(x+e_i)+(1-s)\bar{\phi}_t(x+e_i))_{i\in I}\right)\rmd s.
\end{multline}
Plugging (\ref{tay}) in (\ref{3}), we have
\begin{eqnarray*}
\lefteqn{\frac{\partial}{\partial t}\sum_{y\in\Lambda_\even} (\tilde{\phi_t}(y))^2}\nonumber\\
&=&-2\sum_{x\in\odl}\,\,\,\sum_{\{j\in I,\atop x+e_j\in\Lambda_\even\}}\sum_{k\in I}\tilde{\phi_t}(x+e_k)\tilde{\phi_t}(x+e_j)\int_0^1 D^{j,k}F_x\left((s{\phi}_t(x+e_i)+(1-s)\bar{\phi}_t(x+e_i))_{i\in I}\right)\rmd s\nonumber\\
&=&-2\sum_{x\in\odl}\,\,\,\sum_{\{j\in I,\atop x+e_j\in\Lambda_\even\}}\bigg[(\tilde{\phi_t}(x+e_j))^2\int_0^1 D^{j,j}F_x\left((s{\phi}_t(x+e_i)+(1-s)\bar{\phi}_t(x+e_i))_{i\in I}\right)\rmd s\nonumber\\
&&~~~~~~~~~~~~~~~~~~~+\sum_{k\in I,k\neq j}\tilde{\phi_t}(x+e_k)\tilde{\phi_t}(x+e_j)\int_0^1 D^{j,k}F_x\left((s{\phi}_t(x+e_i)+(1-s)\bar{\phi}_t(x+e_i))_{i\in I}\right)\rmd s\bigg].
\end{eqnarray*}
Using now (\ref{fee0}) for each term $ D^{j,j}F_x\left((s{\phi}_t(x+e_i)+(1-s)\bar{\phi}_t(x+e_i))_{i\in I}\right)$ in the above, we get
\begin{eqnarray}
\label{iexp}
\lefteqn{\frac{\partial}{\partial t}\sum_{y\in\Lambda_\even} (\tilde{\phi_t}(y))^2}\nonumber\\
&=&2\sum_{x\in\odl}\,\,\,\sum_{\{j\in I,\atop x+e_j\in\Lambda_\even\}}\sum_{k\in I,k\neq j}\left[{\tilde{\phi_t}^2(x+e_j)}-\tilde{\phi_t}(x+e_k)\tilde{\phi_t}(x+e_j)\right]\nonumber\\
&&~~~~~~~~~~~~~~~~~~~~~~~~~~~~~\int_0^1 D^{j,k}F_x\left((s{\phi}_t(x+e_i)+(1-s)\bar{\phi}_t(x+e_i))_{i\in I}\right)\rmd s\nonumber\\
&=&2\sum_{x\in\odl}\,\,\,\sum_{\{j,k\in I,j\neq k,\atop x+e_j,x+e_k\in\Lambda_\even\}}\left[{\tilde{\phi_t}^2(x+e_j)}-\tilde{\phi_t}(x+e_k)\tilde{\phi_t}(x+e_j)\right]\nonumber\\
&&~~~~~~~~~~~~~~~~~~~~~~~~~~~~~\int_0^1 D^{j,k}F_x\left((s{\phi}_t(x+e_i)+(1-s)\bar{\phi}_t(x+e_i))_{i\in I}\right)\rmd s\nonumber\\
&&+2\sum_{x\in\odl}\,\,\,\sum_{\{j\in I,\atop x+e_j\in\Lambda_\even\}}\sum_{\{k\in I|\atop x+e_k\in\partial\Lambda_\even\}}\left[{\tilde{\phi_t}^2(x+e_j)}-\tilde{\phi_t}(x+e_k)\tilde{\phi_t}(x+e_j)\right]\nonumber\\
&&~~~~~~~~~~~~~~~~~~~~~~~~~~~~~\int_0^1 D^{j,k}F_x\left((s{\phi}_t(x+e_i)+(1-s)\bar{\phi}_t(x+e_i))_{i\in I}\right)\rmd s,
\end{eqnarray}
where for the second equality we differentiated between $k\in I$ such that $x+e_k\in\Lambda_\even$ and $k\in I$ such that $x+e_k\in\partial\Lambda_\even$. Taking account of the fact that $ D^{j,k}F_x= D^{k,j}F_x$ in the first sum in the last equality above, (\ref{iexp}) becomes
\begin{eqnarray}
\label{iexp1}
\lefteqn{\frac{\partial}{\partial t}\sum_{y\in\Lambda_\even} (\tilde{\phi_t}(y))^2}\nonumber\\
&=&\sum_{x\in\odl}\,\,\sum_{\{j,k\in I, j\neq k|\atop x+e_j,x+e_k\in\Lambda_\even\}}\left[\tilde{\phi_t}(x+e_j)-\tilde{\phi_t}(x+e_k)\right]^2
\int_0^1 D^{j,k}F_x\left((s{\phi}_t(x+e_i)+(1-s)\bar{\phi}_t(x+e_i))_{i\in I}\right)\rmd s \nonumber\\
&&+2\sum_{x\in\odl}\,\,\,\sum_{\{j\in I,\atop x+e_j\in\Lambda_\even\}}\sum_{\{k\in I|\atop x+e_k\in\partial\Lambda_\even\}}\left[{\tilde{\phi_t}^2(x+e_j)}-\tilde{\phi_t}(x+e_k)\tilde{\phi_t}(x+e_j)\right]\nonumber\\
&&~~~~~~~~~~~~~~~~~~~~~~~~~~~~~\int_0^1 D^{j,k}F_x\left((s{\phi}_t(x+e_i)+(1-s)\bar{\phi}_t(x+e_i))_{i\in I}\right)\rmd s\nonumber\\
&\le& -\underline{c}\sum_{b\in\bedl}\left[\nabla\tilde{\phi_t}(b)\right]^2+2 \bar{c}\sum_{b\in\partial\bedl}|\phi_t(y_b)| |\nabla\tilde{\phi_t}(b)|,
\end{eqnarray}
where we used Theorem \ref{rw} and Definition (\ref{rwd}) in the equality in (\ref{iexp1}) to estimate the terms $D^{j,k}F_x\left((s{\phi}_t(x+e_i)+(1-s)\bar{\phi}_t(x+e_i))_{i\in I}\right)$. 
\endpf

~

\noindent \textbf{Coupling Argument}
Suppose that there exist $\mu^\even\in \left(\ext{{\cal S}_\even}\right)_{u}$ and ${\bar{\mu}}^\even\in \left(\ext{{\cal S}_\even}\right)_{v}$ for $u,v\in\RR^d$. For $r>0$, recall the definition of $\chi_{\even,r}$ as given in subsection 3.1. Let us construct two independent $\chi_{\even,r}$-valued 
random variables $\eta_\even=\{\eta_\even(b)\}_{b\in\bed}$ and ${\bar{\eta}}_\even=\{{\bar{\eta}}_\even(b)\}_{b\in\bed}$ on a common probability space $(\Omega,F,P)$ in such a manner that $\eta_\even$ and 
${\bar{\eta}}_\even$ are distributed by $\mu^\even$ and ${\bar{\mu}}^\even$ respectively. We define $\phi_0=\phi^{\eta_\even,0}$ and $\bar{\phi}_0=\phi^{\bar{\eta}_\even,0}$ using the notation in (\ref{19ev}). Let
$\phi_t$ and $\bar{\phi}_t$ be two solutions of the SDE (\ref{sde}) with common Brownian motions having initial data $\phi_0$ and $\bar{\phi}_0$. Let $\eta_{\even,t}$ and
$\bar{\eta}_{\even,t}$ be defined by $\eta_{\even,t}(b):=\nabla\phi(b)$ and $\bar{\eta}_{\even,t}(b):=\nabla\bar{\phi}(b)$, for all $b\in\bed$.
Since $\mu^\even,\bar{\mu}^\even\in {\cal S}_\even$, we conclude that $\eta_{\even,t}$ and
$\bar{\eta}_{\even,t}$ are distributed by $\mu^\even$ and $\bar{\mu}^\even$ respectively, for all $t\ge 0$.

~

\noindent \textbf{Change of Basis}
To adapt the coupling argument from Lemma 2.1 in \cite{FS} to the even bonds, we will use the generator set in $\ed$ outlined below: 
$$
e_{\even,i}=e_{i}+e_{i+1},~i=1,2,\ldots d-1~\mbox{and}~
e_{\even,d}=\left\{
\begin{array}{ll}
e_d-e_1 & d~\mbox{even},\\
e_d+e_1 & d~\mbox{odd}.
\end{array}\right.
$$
Once we have defined this generator set, we can proceed with our arguments. We claim that:
\begin{lemma}
\label{2.1}
There exists a constant $C>0$ independent of $u,v\in\RR^d$ such that
\begin{equation}
\label{0}
{\overline{\lim}}_{T\rightarrow\infty}\frac{1}{T}\int_0^T\sum_{i=1}^{d} \ey^P\left[\left(\eta_{\even,t}(e_{\even,i})-{\bar{\eta}}_{\even,t}(e_{\even,i})\right)^2\right]\rmd t\le C||u-v||^2.
\end{equation}
\end{lemma}
\pf
To prove (\ref{0}), we apply Lemma \ref{lsde} to the differences $\{\tilde{\phi}_t(x):=\phi_t(x)-\bar{\phi}_t(x)\}$ to bound, with the choice $\Lambda_N=[-N,N]^d$, the term
\begin{equation*}
\label{bounderg}
\int_0^T\sum_{x\in\Lambda_N}\ey^P[\tilde{\phi}_t(x)]^2dt.
\end{equation*}
By using shift-invariance in the resulting inequality, we will obtain an upper bound for the term on the left of (\ref{0}). We will next use a special ergodic theorem for co-cycles (see for example Theorem 4 in \cite{co1}), which we can use in our case because $\ed$ is a sub-algebra; we apply it to $\mu^\even\in(\ext {\cal{S_\even}})_{u}$ to obtain
\begin{equation}
\label{ergtheo}
\lim_{\|x\|\rightarrow\infty,\atop x\in\ed}\frac{1}{\|x\|}\|\phi^{\eta_\even,0}(x)-x\cdot u\|_{L^2(\mu^\even)}=0.
\end{equation}
This ergodic theorem will allow us to further estimate the upper bound we have obtained for the term on the left of (\ref{0}), and to obtain the statement of the lemma. The details of the proof, following the same arguments as Lemma 2.1 from \cite{FS}, will be omitted and are left to an interested reader.
\endpf

\subsection{Step 2: Uniqueness of ergodic component for $(\Z^d)^*$}
\textbf{Proof of Theorem \ref{ergall}}
Let $u\in\RR^d$. Suppose now that there exist $\mu,\bar{\mu}\in {\cal G}_\beta(H)$ ergodic and shift-invariant such that $\ey_{\mu}(\eta(b))=\ey_{\bar{\mu}}(\eta(b))=\langle u, y_b-x_b\rangle$ for all bonds $b=(x_b,y_b)\in\bzd$.   Note now that $\ey_{\mu^\even}(\eta_\even(b))=\ey_{\bar{\mu^\even}}(\eta_\even(b))=\langle u, y_b-x_b\rangle$ for all bonds $b=(x_b,y_b)\in\bed$.

From Lemma~\ref{muodd} and with the same notation as there, we get that $\mu^\even,\bar{\mu}^\even\in {\cal G}_\even(H^{\even})$. As for all $\eta_\even\in\chi_\even$, with $\eta_\even(b)=\phi(y_b)-\phi(x_b)$, $b=(x_b,y_b)\in\bed$, we can write 
$\eta_\even(b)=\eta(b_1)+\eta(b_2)$, $b_1,b_2\in\bzd$,
shift-invariance and ergodicity under the even shifts for $\mu^\even,\bar{\mu}^\even$ follow immediately from the similar properties for $\mu,\bar\mu$. Therefore $\mu^\even,{\bar{\mu}}^\even\in\left(\ext{{\cal S}_\even}\right)_{u}$, so we can apply Theorem~\ref{ergeven1} to get $\mu^\even=\bar{\mu}^\even$.
Then for any $A\in {{\cal F}}_\bzd$, we have from Lemma~\ref{led} that $\ey_{\mu}(1_A|{{\cal F}}_\bed)=\ey_{\bar\mu}(1_A|{{\cal F}}_\bed)$ and we have
$$\mu(A)=\ey_{\mu}(1_A)=\ey_{\mu}(\ey_{\mu}(1_A|{{\cal F}}_\bed))=\ey_{\bar{\mu}}(\ey_{\mu}(1_A|{{\cal F}}_\bed))=\ey_{\bar{\mu}}(\ey_{\bar{\mu}}(1_A|{{\cal F}}_\bed))=\ey_{\bar{\mu}}(A)=\bar{\mu}(A).$$
\endpf

\subsection{Existence of ergodic component on $(\Z^d)^*$}

Tightness of the family $\{\mu_{\Lambda,\xi}\}_{\Lambda\subset\zd,\xi\in\chi}$ is known for strictly convex potentials with quadratic growth at $\infty$ (see for example Section 4.4 in \cite{FSL}). Therefore a limiting measure exists in this case by taking $|\Lambda|\rightarrow\infty$ along a suitable sub-sequence. For non-convex potentials satisfying (\ref{tag2}) and such that $U''(s)\le C_2$ for all $s\in\RR$, tightness of the family $\{\mu_{\Lambda,\xi}\}_{\Lambda\subset\zd,\xi\in\chi}$ and existence of the limiting measure are shown in \cite{CK} in a more general situation (see Lemmas 3.6 and 3.7 and Proposition 3.8 from \cite{CK}). 

To automatically ensure shift invariance, we will construct below shift-invariant Gibbs measures through the use of \textit{periodic boundary conditions}. For this reason, take $N\in\N$ and let $\TT_N^d=(\ZZ/N\ZZ)^d$ be the lattice torus in $\zd$. As before, $(\TT_N^d)^*$ denotes the set of directed bonds in $\TT_N^d$ and $\chi_{\TT_N^d}$ denotes the set of all $\eta\in\RR^{(\TT_N^d)^*}$ which satisfy the plaquette condition.

\begin{lemma}
\label{exist}
Let $U$ be such that it satisfies (\ref{tag2}) and such that $U''(s)\le C_2$ for all $s\in\RR$. Then for every $u\in\RR^d$ there exists at least one shift-invariant $\mu_u\in {\cal G}_\beta(H)$ with a given tilt $u\in\RR^d$.
\end{lemma}
\pf
For the proof of existence of shift-invariant $\nabla\phi$-Gibbs measures we proceed as in the proof of Theorem 3.2 from \cite{FS}. To avoid that only the state with tilt $u=0$ could be
constructed, we note that boundary conditions with definite
tilt $u$ are identical to boundary conditions $u=0$ but with the shifted potential $U (\cdot+u_i)$
for a bond directed along $e_i,i\in I$.

Fix $u\in\RR^d$ and let
\begin{equation}
\tilde\mu_{N,u}( d\tilde \eta)=\frac{1}{Z_{N,u}} 
\exp\bigl( 
-\beta\sum_{b \in (\TT_N^d)^*}U(\tilde \eta(b)+u_b)\bigr)\rmd\tilde \eta_N.
\end{equation}
Here $\rmd\tilde \eta_N$ is the uniform measure on the affine space $\chi_{\TT_N^d}$, $Z_{N,u}$ is the normalization and $u_b:=\pm u_i$ for $b=(x\pm e_i,x),x\in\zd,i\in \{1,\ldots, d\}$. The law of $\{\eta(b):=\tilde\eta(b)+u_b\}$ under $\tilde\mu_{N,u}$ is denoted by $\mu_{N,u}$.

Consider 
\begin{equation}
\label{eq1}
\limsup_{N\uparrow\infty}\frac{1}{|\TT_N^d|}\log 
\tilde\mu_{N,u}\Bigl(
\exp\Bigl(\gamma\sum_{b \in (\TT_N^d)^*}(\tilde \eta(b))^2
\Bigr),
\end{equation}
where $\gamma>0$ will be chosen later. We will find next an upper bound for this expression.
\begin{equation*}
\label{eq2}
\tilde\mu_{N,u}\Bigl(
\exp\Bigl(
\gamma\sum_{b \in (\TT_N^d)^*}(\tilde \eta(b))^2
\Bigr)
\Bigr)=\frac{\int \exp\bigl( 
-\beta\sum_{b \in (\TT_N^d)^*}U(\tilde \eta(b)+u_b)+\gamma\sum_{b \in (\TT_N^d)^*}(\tilde \eta(b))^2
\bigr) d\tilde \eta_N}{
\int \exp\bigl( 
-\beta\sum_{b \in (\TT_N^d)^*}U(\tilde \eta(b)+u_b)
\bigr) d\tilde \eta_N}.
\end{equation*}
Using the assumption on the potential, 
$U(s)\leq C_2 s^2+U(0)$ and $U(s)\geq A s^2 -B$, 
this expression is bounded from above 
by
$$e^{\beta(B-U(0)) |\TT_N^d|}\frac{\int \exp\bigl( 
-\beta\sum_{b \in \TT_N^d}A(\tilde \eta(b)+u_b)^2 +\gamma\sum_{b \in \TT_N^d}(\tilde \eta(b))^2
\bigr) d\tilde \eta_N}{
\int \exp\bigl( 
-\beta\sum_{b \in (\TT_N^d)^*}C_2(\tilde \eta(b)+u_b)^2\bigr) d\tilde \eta_N}.$$
By Remark~\ref{equivzd1}, we can express the uniform integration over gradient fields as an integration over the 
fields $\tilde\phi(x)=\phi(x)-u\cdot x$, and the above expression is equal to
\begin{eqnarray}
\label{eq3}
\lefteqn{e^{(\beta B-\beta U(0)) |\TT_N^d|}\frac{1}{
\int \exp\bigl( 
-\beta C_2
\sum_{x \in \TT_N^d\atop i\in I}(\tilde \phi(x)-\tilde\phi(x+e_i)+u_i )^2\bigr) d\tilde \phi_{\TT_N^d\setminus \{0\}}\delta_0(d\tilde\phi(0))}}\nonumber\\
&&\times \int \exp\bigl( 
-A\beta\sum_{x \in \TT_N^d\atop i\in I}(\tilde \phi(x)-\tilde\phi(x+e_i)+u_i )^2+\gamma\sum_{x \in \TT_N^d\atop i\in I}(\tilde \phi(x)-\tilde\phi(x+e_i))^2\bigr) d\tilde \phi_{\TT_N^d\setminus \{0\}}\delta_0(d\tilde\phi(0)).\nonumber\\
\end{eqnarray}
But
\begin{eqnarray}
\label{eq4}
\lefteqn{-A\beta\sum_{x \in \TT_N^d\atop i\in I}\bigl((\tilde \phi(x)-\tilde\phi(x+e_i)+u_i )^2+ \gamma\sum_{x \in \TT_N^d\atop i\in I}(\tilde \phi(x)-\tilde\phi(x+e_i))^2\bigr)}\nonumber\\
&=&-(A\beta-\gamma)\sum_{x \in \TT_N^d\atop i\in I}(\tilde \phi(x)-\tilde\phi(x+e_i))^2-A\beta|\TT_N^d|\sum_{i\in I}u_i ^2.
\end{eqnarray}
Let $\gamma<A\beta$ be arbitrarily fixed. Plugging (\ref{eq3}) and (\ref{eq4}) in (\ref{eq1}) and integrating out, we obtain for some $C(\beta,A,C_2,u)>0$
$$\limsup_{N\uparrow\infty}\frac{1}{|\TT_N^d|}\log 
\tilde\mu_{N,u}\Bigl(
\exp\Bigl(\gamma\sum_{b \in (\TT_N^d)^*}(\tilde \eta(b))^2
\Bigr)<C(\beta,A,C_2,u)<\infty. $$
In particular, due to the shift-invariance of the family $(\tilde\mu_{N,u})_{N\in\N}$ on $\TT_N^d$, we get from the above for all bonds $b$
$$\limsup_{N\uparrow\infty}\tilde\mu_{N,u}((\tilde \eta(b))^2)<C(\beta,A,C_2,u)<\infty, $$
which implies tightness of the family $(\tilde\mu_{N,u})_{N\in\N}$.
\endpf

\begin{thm} \textbf{(Existence of ergodic component on $(\Z^d)^*$)}
\label{existerg}
Let $U=V+g$, where $U$ satisfy (\ref{tag2}) and $V$ and $g$ satisfy (\ref{vc}) and  (\ref{tag5}). Then for every $u\in\RR^d$, there exists at least one ergodic, shift-invariant $\mu_u\in {\cal G}_\beta(H)$ with a given tilt $u\in\RR^d$.
\end{thm}
\pf
Existence of shift-invariant $\mu\in P_2({\chi})$ with given tilt $u\in\RR^d$ is assured for our non-convex class of potentials by Lemma \ref{exist}; nevertheless, existence of an ergodic and shift-invariant $\mu_u\in P_2({\chi})$ with given tilt $u\in\RR^d$ is not assured for non-convex potentials. However, due to the strict convexity of the $F_x$ potentials, we can use the Brascamp-Lieb inequality and a similar reasoning to the one of Theorem 3.2 in \cite{FS}, to easily show the existence, for every $u\in\RR^d$, of \textbf{at least one} $\mu_u\in {\cal G}_\beta(H)$ ergodic and shift-invariant and with tilt $u\in\RR^d$.
\endpf

\section{Decay of Covariances}
In this section,
we extend the covariance estimates of  \cite{DD} to the class of non-convex potentials $U=V+g$ which satisfy (\ref{tag2}) such $V$ 
and $g$ satisfy (\ref{vc}) and (\ref{tag5}).

Recall that $F\in C^1_b(\chi_r)$, where $C^1_b(\chi_r)$ denotes the set of differentiable functions depending on finitely many 
coordinates with bounded derivatives and where $\chi_r$ was defined in subsection 1.2.2. Using now $\eta,\eta'\in\chi_\even$ in (\ref{partialbond}), we define $\partial_{b_\even} F$ and $||\partial_{b_\even}F||_{\infty}$ similarly for $b_\even\in\bed$ as we did for $b\in\bzd$. Before proving Theorem \ref{cov}, we make a remark which we will use in our proof.
\begin{rem}\normalfont
Take $b_\even=(x+e_l,x+e_j)\in\bed$. In view of the definition, we have 
\begin{equation}
\label{betob} 
||\partial_{b_\even}F||_{\infty}=\sup_{\eta\in\chi_\even}\left|\partial_{b_\even} F(\eta)\right|\le\sum_{b\in\bzd: b\sim b_\even}\sup_{\eta\in\chi}\left|\partial_bF(\eta)\right|=\sum_{b\in\bzd:  b\sim b_\even}||\partial_bF||_{\infty},
\end{equation}
where $b\sim b_\even$ are those $b=(x,x+e_s)\in (\zd)^*,x\in\od$, such that $s\in\{l,j\}$. 
\end{rem}
\textbf{Proof of Theorem \ref{cov}}
We have 
\begin{eqnarray}
\label{alt1}
\cov_{\mu_u}(F(\eta),G(\eta))&=&\ey_{\mu_u}\left[\cov_{\mu_u}(F(\eta),G(\eta)|{\cal F}_\bed)\right]\nonumber\\
&&+\cov_{\mu_u}\left(\ey_{\mu_u}[F(\eta)|{\cal F}_\bed],\ey_{\mu_u}[G(\eta)|{\cal F}_\bed]\right),
\end{eqnarray}
where by Corollary \ref{tildemuaa} and with the same notations, we have for a fixed $k\in I$
$$\ey_{\mu_u}\left(F|{\cal F}_{\bed}\right)(\eta)=\int F\left(\left(\eta_\even(b)-\phi(x)\right)_{b\in {\mathcal B}(x,k),x\in\od}\right)\prod_{x\in\od}\mu_{x,\eta_\even}^k(\rmd\phi(x));$$
a similar formula holds for $G$. Note that under $\mu_u(\,\cdot\,|{\cal F}_\bed)$, the gradient vectors $((\nabla\phi_i(x))_{i\in I})_{ x\in\od}$ are independent for all $x\in\od$. In view of this and of the above formula, under $\mu_u(\,\cdot\,|{\cal F}_\bed)$ the gradients $(\nabla\phi_i(x),i\in I, x\in\od)$ are pairwise positive quadrant dependent. That means that for all $x,y\in\od, i,j\in I$, with either $x\neq y$ or $i\neq j$, we have 
\begin{multline}
\label{65}
\ey_{\mu_u}\left(1_{\left(\nabla\phi_i(x)>a_i,\nabla\phi_j(y)>a_j\right)}|{\cal F}_{\bed}\right)(\eta)\\
\ge\ey_{\mu_u}\left(1_{\left(\nabla\phi_i(x)>a_i\right)}|{\cal F}_{\bed}\right)(\eta)\ey_{\mu_u}\left(1_{\left(\nabla\phi_j(y)>a_j\right)}|{\cal F}_{\bed}\right)(\eta),~\forall\eta\in\chi~\mbox{ and}~ \forall a_i,a_j\in\RR.
\end{multline}
To show this, note first that the inequality is true with equal sign for all $x,y\in\od, i,j\in I,x\neq y$, due to the independence of the gradient vectors. For the case with $x=y\in\od,i,j\in I,$ the left-hand side of (\ref{65}) becomes in view of Lemma \ref{tildemuaa}
\begin{eqnarray*}
\lefteqn{\ey_{\mu_u}\left(1_{\left(\nabla\phi_i(x)>a_i,\nabla\phi_j(y)>a_j\right)}|{\cal F}_{\bed}\right)(\eta)}\\
&=&\int 1_{\left(\phi(x+e_i)-\phi(x+e_k)-\phi(x)>a_i,\phi(x+e_j)-\phi(x+e_k)-\phi(x)>a_j\right) }(\phi(x))\prod_{x\in\od}\mu_{x,\eta_\even}^k(\rmd\phi(x))\\
&=&\int 1_{\left(\phi(x)<\min\{\phi(x+e_i)-\phi(x+e_k)-a_i,\phi(x+e_j)-\phi(x+e_k)-a_j\}\right) }(\phi(x))\prod_{x\in\od}\mu_{x,\eta_\even}^k(\rmd\phi(x))\\
&=&\min\bigg(\int 1_{\left(\phi(x)<\phi(x+e_i)-\phi(x+e_k)-a_i\right) }(\phi(x))\prod_{x\in\od}\mu_{x,\eta_\even}^k(\rmd\phi(x)),\\
&&\,\,\,\,\,\,\,\,\,\,\,\,\,\,\,\,\,\,\int 1_{\left(\phi(x)<\phi(x+e_j)-\phi(x+e_k)-a_j\right) }(\phi(x))\prod_{x\in\od}\mu_{x,\eta_\even}^k(\rmd\phi(x)) \bigg)\\
&=&\min\left(\ey_{\mu_u}\left(1_{\left(\nabla_i\phi(x)>a_i\right)}|{\cal F}_{\bed}\right)(\eta),\ey_{\mu_u}\left(1_{\left(\nabla_j\phi(x)>a_j\right)}|{\cal F}_{\bed}\right)(\eta)\right)\\
&\ge&\ey_{\mu_u}\left(1_{\left(\nabla\phi_i(x)>a_i\right)}|{\cal F}_{\bed}\right)(\eta)\ey_{\mu_u}\left(1_{\left(\nabla\phi_j(y)>a_j\right)}|{\cal F}_{\bed}\right)(\eta), 
\end{eqnarray*}
so the inequality holds. Note now that Lemma~3.1 from \cite{JD} can be adapted to the case with pairwise positive quadrant dependent random variables. The reason for this is that the main ingredient used in Lemma 3.1, Rosenthal's inequality, holds for the case with pairwise positive quadrant dependent random variables (see, for example, Corollary 1 from \cite{ros} for a statement of Rosenthal's theorem in this case). Given (\ref{alt1}), the rest of the 
argument from Lemma 3.1 can be easily adapted to our case; therefore, there exists $c>0$ such that 
\begin{eqnarray}
\label{cv11b}
\left|\cov_{\mu_u}(F(\eta),G(\eta)|{\cal F}_\bed)\right|&\le& c\sum_{b\in\bzd}||\partial_b F||_{\infty}||\partial_b G||_{\infty}\var_{\mu_u}(\nabla\phi(b)|{\cal F}_\bed)\nonumber\\
&\le&c'\tau^2\sum_{b\in\bzd}||\partial_b F||_{\infty}||\partial_b G||_{\infty},
\end{eqnarray}
where the first inequality is an application of the adaptation of Lemma~3.1 in \cite{JD}, and for the second inequality we used (\ref{newvga}).
Note that, due to the fact that the random walk representation holds, Theorem~6.2 from \cite{DD} can be adapted to the case of the infinite even lattice with strictly convex potential; thus, a decay of covariance statement, similar to the one in Theorem \ref{cov}, holds for the even setting. In view of Lemma~\ref{muodd}, there exists $c''>0$ such that
\begin{equation}
\label{hatcv1}
\left|\cov_{\mu_u}(\hat{F},\hat{G})\right|\le c''\sum_{b_\even,b'_\even\in\bed}\frac{||\partial_{b_\even}\hat{F}||_{\infty}||\partial_{b'_\even}\hat{G}||_{\infty}}{1+\|x_{\even}-x'_{\even}\|^d},
\end{equation}
where $\hat{F}=\ey_{\mu_u}[F(\eta)|{\cal F}_\bed]$ and $\hat{G}=\ey_{\mu_u}[G(\eta)|{\cal F}_\bed]$. We need to estimate now $\partial_{b_\even}\hat{F}$ and $\partial_{b_\even}\hat{G}$. But
\begin{multline}
\label{bonder}
\partial_{b_\even}\hat{F}=\partial_{b_\even}\ey_{\mu_u}[F(\eta)|{\cal F}_\bed]=\ey_{\mu_u}[\partial_{b_\even} F(\eta)|{\cal F}_\bed]\\
-\cov_{\mu_u}\left(F(\eta),\partial_{b_\even}\left(\sum_{x\in\od}\sum_{b\in {\mathcal B}(x,k)} U\left(\eta_\even(b)-\phi(x)\right)\right)\bigg|{\cal F}_\bed\right),
\end{multline}
from which, by using also (\ref{betob})
\begin{equation}
\label{hatcv2}
|\partial_{b_\even}\hat{F}|\le \sum_{b:b\sim b_\even}||\partial_b F||_{\infty}+\bigg|\cov_{\mu_u}\big(F(\eta),\sum_{x\in\od,\atop b_\even\in {\mathcal B}(x,k)} U'\left(\eta_\even(b_\even)-\phi(x) \right)\bigg|{\cal F}_\bed\big)\bigg|.
\end{equation}
Applying (\ref{cv11b}) to the covariance in (\ref{hatcv2}) and using $|U''|\le C_0+C_2$ and (\ref{newvga}), we get for some $c'''>0$
\begin{multline}
\label{cv11bb}
\left|\cov_{\mu_u}\left(F(\nabla\phi),\partial_{b_\even}\left(\sum_{x\in\od}\sum_{b\in {\mathcal B}(x,k)} U\left(\eta_\even(b)-\phi(x)\right)\right)\bigg|{\cal F}_\bed\right)\right|\\
\le 2d c'''(C_0+C_2)||\partial_{b_\even} F||_{\infty}\var_{\mu_u}(\eta(b)|{\cal F}_\bed)
\le\tilde{c}||\partial_{b_\even} F||_{\infty}.
\end{multline}
The statement of the theorem follows now from (\ref{hatcv2}), (\ref{cv11bb}), (\ref{cv11b}), (\ref{hatcv1}) and (\ref{betob}).
\endpf

\section{Central Limit Theorem}
We will extend next in Theorem \ref{clt} the scaling limit results from \cite{gos} to our class of potentials.
\\
\textbf{Proof of Theorem \ref{clt}}
It suffices to prove that for all $i\in I$
$$S_{\epsilon,i}(f)=\epsilon^{d/2}\sum_{x\in\zd}f(x\epsilon)(\nabla_i\phi(x) -u_i) \Rightarrow N(0,\sigma^2_{u,i}(f))~~\mbox{as}~~\epsilon\rightarrow 0.$$
Note that
\begin{eqnarray*}
S_{\epsilon,i}(f)&=&\epsilon^{d/2}\sum_{x\in\zd}f(x\epsilon)\left[\phi(x+e_i)-\phi(x) -u_i\right]=\epsilon^{d/2}\sum_{x\in\ed}f(x\epsilon)\left[\phi(x+2e_i)-\phi(x) -2u_i\right]\nonumber\\
&&-\epsilon^{d/2}\sum_{x\in\ed}f(x\epsilon)\left[\phi(x+2e_i)-\phi(x+e_i) -u_i\right]+\epsilon^{d/2}\sum_{x\in\od}f(x\epsilon)\left[\phi(x+e_i)-\phi(x) -u_i\right]\nonumber\\
&=&\epsilon^{d/2}\sum_{x\in\ed}f(x\epsilon)\left[\phi(x+2e_i)-\phi(x) -2u_i\right]\nonumber\\
&&+\epsilon^{d/2}\sum_{x\in\ed}\Big[f((x+e_i)\epsilon)-f(x\epsilon)\Big]\left[\phi(x+2e_i)-\phi(x+e_i) -u_i\right]= S^e_{\epsilon}(f)+R_{\epsilon}(f).
\end{eqnarray*}
We can show the CLT for $S^e_{\epsilon,i}(f)$ since the summation is concentrated on the even sites; the proof uses the same arguments as in \cite{gos} 
and is based on the random walk representation, as explained in Remark \ref{rwrepres}. Also, since by Theorem~\ref{cov}
$$\big|\cov_{\mu_u}(\nabla_i\phi(x),\nabla_j\phi(y))\big|\le \frac{C}{(\|x-y\|+1)^d},$$
we have
\begin{eqnarray*}
\var_{\mu_u}(R_{\epsilon,i}(f))&\le&\epsilon^{d} \sum_{x,y\in\ed} |\nabla_i f(x\epsilon)| |\nabla_i f(y\epsilon)||\cov_{\mu_u}(\phi(x+e_i)-\phi(x),\phi(y+e_i)-\phi(y))\big|\nonumber\\
&\le&\epsilon^{d} \sum_{x,y\in\ed} |\nabla_i f(x\epsilon)| |\nabla_i f(y\epsilon)|
\frac{C}{(\|x-y\|+1)^d},
\end{eqnarray*}
where $\nabla_i f(x\epsilon)=f((x+e_i)\epsilon)-f(x\epsilon)$.
Expanding $f((x+e_i)\epsilon)$ around $x\epsilon$ by the Mean Value Theorem, we have $\nabla_i f(x\epsilon)=D^i f(a)\epsilon,$
for some $a\in\RR^d$. As $f\in C_0^{\infty}(\RR^d)$, there exist $M,N>0$ such that for all $x\in\RR^d$ with $|\epsilon x|\le N$ we 
have $f(\epsilon x)\le M$, $|D^i f(\epsilon x)|\le M$ and both functions equal to $0$ for $|\epsilon x|>N$. Therefore
\begin{eqnarray*}
\var_{\mu_u}(R_{\epsilon,i}(f))&\le&\sum_{x,y\in\ed,\atop |\epsilon x|\le N, |\epsilon y|\le N}\frac{\epsilon^{d+2}M^2 C}{(\|x-y\|+1)^d}\le\epsilon^{d+2}M^2 C\sum_{y\in\ed,\atop |\epsilon y|\le N}\int_{-\frac{N}{\epsilon}}^{\frac{N}{\epsilon}}\ldots\int_{-\frac{N}{\epsilon}}^{\frac{N}{\epsilon}}\frac{\rmd x_1\rmd x_2\ldots\rmd x_d}{\left(\sum_{i=1}^d|x_i-y_i|+1\right)^d}\nonumber\\
&\le&\epsilon^2 C(d,N,M)\log\left(1+2dN/\epsilon\right)\le 2dN C(d,N,M)\epsilon,
\end{eqnarray*}
where $C(d,N,M)$ is a positive constant depending on $d,M$ and $N$.
It follows that $R_{\epsilon,i}(f)\to 0$ in probability as $\epsilon\to 0$.
\endpf

\section{Surface tension}

We will extend here the surface tension strict convexity results from \cite{FS} and \cite{DGI} to the family of non-convex potentials satisfying (\ref{tag2}), (\ref{vc}) and (\ref{tag5}). 

Take $N\in\N$ and let $\TT_N^d=(\ZZ/N\ZZ)^d$ be the lattice torus in $\zd$ and let $u\in\RR^d$. Then, we define the surface tension on the torus $\TT_N^d$ as
$$\sigma_{\TT_N^d}^{\beta}(u)=-\frac{1}{|\TT_N^d|}\log\frac{Z_{\TT_N^d}^{\beta}(u)}{Z_{\TT_N^d}^{\beta}(0)},~\mbox{with}~Z_{\TT_N^d}^{\beta}(u)=\int_{\RR^{\TT_N^d}} 
\exp(-\beta H_{\TT_N^d}(\phi,u))\prod_{x\in 
\TT_N^d\setminus\{0\}} \rmd\phi(x)$$
and where $H_{\TT_N^d}$ is given by 
$$H_{\TT_N^d}(\phi,u)=\sum_{x\in\TT_N^d}\sum_{i=1}^dU(\btd_{i}\phi(x)+u_i) =\sum_{x\in\TT_N^d}\sum_{i=1}^d\left[V(\btd_{i}\phi(x)+u_i)+g(\btd_{i}\phi(x)+u_i)\right].$$
We define $u_{-i}=-u_i$ for $i=1,2,\ldots, d$.
Take now $N$ to be even. Just as in the previous sections, let us label the vertices of the torus as odd and even; let the set of odd vertices on the torus be $\tnod$ and the set of even vertices be $\tnev$. Then we can of course first integrate all the odd coordinates and:
\begin{eqnarray*}
Z_{\TT_N^d}^{\beta}(u)&=&\int_{\RR^{\edn}}\left(\int_{\RR^{\tnod}}\exp(-\beta H_{\TT_N^d}(\phi,u)\prod_{x\in\tnod}\rmd\phi(x)\right)\prod_{x\in\tnev\setminus\{0\}}\rmd\phi(x)\nonumber\\
&=&
\int_{\RR^{\tnev}}\exp(-\beta H^\even_{\tnev}(\phi,u))\prod_{x\in\tnev\setminus\{0\}}d\phi(x),
\end{eqnarray*}
where, similarly to (\ref{eqnW})
$$H_{\EEN}^\even(\phi,u)=\sum_{x\in\tnod}F_x((\phi(x+e_i))_{i\in I},u),~I=\{-d,\ldots,d\}\setminus\{0\},$$
with
$$F_x((\phi(x+e_i))_{i\in I},u)=-\log\int_{\RR} e^{-\beta\sum_{i\in I} U(\btd_{i}\phi(x)+u_i)}\rmd\phi(x).$$
Then, defining the \textbf{even} surface tension on $\tnev$ as
$$\sigma_{\EEN}^{\beta}(u)=-\frac{1}{|\tnev|}\log\frac{Z_{\EEN}^{\beta}(u)}{Z_{\EEN}^{\beta}(0)},~\mbox{with}~ Z_{\EEN}^{\beta}(u)=\int_{\RR^{\tnev}} 
\exp(-\beta H_{\EEN}^\even(\phi,u))\prod_{x\in 
\tnev\setminus\{0\}} \rmd\phi(x),$$
we obtain the following result by integrating out the odds
\begin{lemma}
\label{evtor}
$$\sigma_{\EEN}^{\beta}(u)=\frac{1}{2}\sigma_{\TT_N^d}^{\beta}(u).$$
\end{lemma}
We will next prove strict convexity for the even surface tension, uniformly in $N$ even.
\begin{thm} [Strict convexity of the even surface tension]
\label{evconv} 
Suppose that $V,g\in C^2(\RR)$ such that they satisfy (\ref{tag2}), (\ref{vc}) and (\ref{tag5}). Then, for all $N=2k$, we have
\begin{equation}
D^2\sigma_{\TT_N^d}^{\beta}(u)=2 D^2\sigma_{\EEN}^{\beta}(u)\ge 4d\beta^2c_l Id,~\forall~u\in\RR^d,
\end{equation}
where $c_l$ is given in (\ref{loboundcov}). That is, the even surface tension is uniformly strictly convex in $u\in\RR^d$, uniformly in all $N$ even.
\end{thm}
\pf
Since $H^{\even}$ fulfills the random walk representation condition by Theorem~\ref{rw}, $F_x$ are uniformly convex and we can apply Lemma 3.2 in \cite{CDM} to $\sigma_{\EEN}^{\beta}(u)$, to get the statement of our theorem.
\endpf
Note now that by the same reasoning as in \cite{FS}, we can prove the existence of
$$\sigma^{\beta}(u)=\lim_{|\TT_N^d|\rightarrow\infty}\sigma_{\TT_N^d}^{\beta}(u).$$ 
Together with Theorem \ref{evconv}, this gives
\begin{thm} [Strict convexity of the surface tension]
\label{evconv1} 
Suppose that $V,g\in C^2(\RR)$ such that they satisfy (\ref{tag2}), (\ref{vc}) and (\ref{tag5}). Then the surface tension $\sigma^{\beta}(u)$ is strictly convex in $u\in\RR^d$.
\end{thm}

\section{Appendix}

Due to the fact that Example 3.2 (a) has been the subject of two other papers in the area (see \cite{BK} and \cite{BS}), we will provide here a sketch of the explicit computations for this example, which provide us with the $\frac{p}{1-p}<O\left(\left(\frac{k_2}{k_1}\right)^{1/2}\right)$ order. The explicit computations are worth separate consideration, as they don't follow from Theorem \ref{rw}. As before, it is sufficient to estimate $\cov_{\nu_{x,\psi}}\left(U'(\nabla_{i}\phi(x)),U'(\nabla_{j}\phi(x))\right)$, for all $x\in\od$ and $i,j\in I,i\neq j$.  

Denote by $\theta_k:=\phi(x+e_k),k=1,\ldots,4$. Let $\Xi:=\{(\alpha,{\bar{\alpha}})|\alpha=(\alpha_1,\ldots,\alpha_4), \bar{\alpha}=(1-\alpha_1,\ldots,1-\alpha_4)\}, ~\mbox{with}~\alpha_k\in\{0,1\}, k=1,\ldots, 4\}$.

Since $U\ge\bar{c}~k_2 $ outside of a domain $[-\frac{{\tilde{c}}}{\sqrt{k_1-k_2}},\frac{{\tilde{c}}}{\sqrt{k_1-k_2}}]$, for some ${\tilde{c}}>0$ and for some $\bar{c}>0$, we take $V,g$ to be defined as in (\ref{3.2a}) on $[-\frac{{\tilde{c}}}{\sqrt{k_1-k_2}} ,\frac{{\tilde{c}}}{\sqrt{k_1-k_2}}]$ and $V:=U, g:=0$, on the complement set.
By the same reasoning as in (\ref{vi}), (\ref{gi}) and (\ref{vx3a}) from Theorem \ref{rw}, we know that the terms $\cov_{\nu_{x,\psi}}(V',V')$ and $\cov_{\nu_{x,\psi}}(g',g')$  are positive terms, while the terms $\cov_{\nu_{x,\psi}}(V',g')$ are negative terms. Using the same reasoning as in example 3.2 (a), we get that 
\begin{equation}
\label{lowboundbk}
\cov_{\nu_{x,\psi}}\left(V'(\nabla_{i}\phi(x)),V'(\nabla_{j}\phi(x))\right)\ge \bar{c}~k_2.
\end{equation}
We will next try to bound from below the negative part of $\cov_{\nu_{x,\psi}}\left(U'(\nabla_{i}\phi(x)),U'(\nabla_{j}\phi(x))\right)$. Note first that, by a reasoning similar to (\ref{vx2}), we get for the negative part
\begin{multline}
\label{vx22}
\cov_{\nu_{x,\psi}}\left(g'(\nabla_j\phi(x)),V'(\nabla_{i}\phi(x))\right)\ge\cov_{\nu_{x,\psi}}\left(g'(\nabla_j\phi(x)),\sum_{k\in I} V'(\nabla_{k}\phi(x))\right)\\
=\frac{1}{2}\ey_{\nu_{x,\psi}}\left(g''(\nabla_{j}\phi(x))\right)-\cov_{\nu_{x,\psi}}\left(g'(\nabla_{j}\phi(x)),\sum_{k\in I}g'(\nabla_{k}\phi(x))\right)\ge\frac{1}{2}\ey_{\nu_{x,\psi}}\left(g''(\nabla_{j}\phi(x))\right).
\end{multline}
We next bound $\ey_{\nu_{x,\psi}}\left(-g''(\nabla_{j}\phi(x))\right)$ from above, where by (\ref{3.2a}) 
$$g''(s)=-\frac{p(1-p)(k_1-k_2)^2s^2}{p^2e^{-(k_1-k_2)\frac{s^2}{2}}+2p(1-p)+(1-p)^2 e^{(k_1-k_2)\frac{s^2}{2}}}\ge -\frac{p(k_1-k_2)^2s^2}{(1-p) e^{(k_1-k_2)\frac{s^2}{2}}}$$
on $[-\frac{{\tilde{c}}}{\sqrt{k_1-k_2}},\frac{{\tilde{c}}}{\sqrt{k_1-k_2}}]$ and $0$ otherwise. Therefore
$$\ey_{\nu_{x,\psi}}\left(-g''(\nabla_{j}\phi(x))\right)\le\frac{p}{1-p}(k_1-k_2)^2\frac{\int_{-\frac{{\tilde{c}}}{\sqrt{k_1-k_2}}}^{\frac{{\tilde{c}}}{\sqrt{k_1-k_2}}}(s-\theta_j)^2e^{-(k_1-k_2)\frac{s^2}{2}} e^{-\sum_{k=1}^4 U(s-\theta_k)}\rmd s}{\int_{\RR} e^{-\sum_{k=1}^4 U(s-\theta_k)}\rmd s},$$
where $U(s)=-\log\left(pe^{-k_1\frac{s^2}{2}}+(1-p)e^{-k_2\frac{s^2}{2}}\right)$. Then
\begin{eqnarray}
\label{40}
\lefteqn{\ey_{\nu_{x,\psi}}\left(-g''(\nabla_{j}\phi(x))\right)}\nonumber\\
&\le&\frac{p}{1-p}(k_1-k_2)^2\frac{\int_{-\frac{{\tilde{c}}}{\sqrt{k_1-k_2}}}^{\frac{{\tilde{c}}}{\sqrt{k_1-k_2}}}(s-\theta_j)^2e^{-(k_1-k_2)\frac{s^2}{2}} \prod_{k=1}^4 \left(p e^{-k_1 \frac{(s-\theta_k)^2}{2}}+(1-p)e^{-k_2\frac{(s-\theta_k)^2}{2}}\right)\rmd s}{\int\prod_{k=1}^4\left(p e^{-k_1 \frac{(s-\theta_k)^2}{2}}+(1-p)e^{-k_2\frac{(s-\theta_k)^2}{2}}\right)\rmd s}\nonumber\\
&=&\frac{p}{1-p}(k_1-k_2)^2\frac{\sum_{(\alpha,{\bar{\alpha}})\in\Xi} \int_{-\frac{{\tilde{c}}}{\sqrt{k_1-k_2}}}^{\frac{{\tilde{c}}}{\sqrt{k_1-k_2}}}(s-\theta_j)^2e^{-(k_1-k_2)\frac{s^2}{2}}I(k_1,k_2,\alpha,{\bar{\alpha}})\rmd s}{\sum_{(\alpha,{\bar{\alpha}})\in\Xi}\int I(k_1,k_2,\alpha,{\bar{\alpha}})\rmd s},
\end{eqnarray}
where $I(k_1,k_2,\alpha,{\bar{\alpha}}):=p^{\sum_{k=1}^4\alpha_k}(1-p)^{\sum_{k=1}^4{\bar{\alpha}}_k} e^{-k_1\sum_{k=1}^4\alpha_k \frac{(s-\theta_k)^2}{2}-k_2\sum_{k=1}^4{\bar{\alpha}}_k \frac{(s-\theta_k)^2}{2} }$, and where (\ref{40}) is a sum of sixteen Gaussian integrals. Define for $(\alpha,{\bar{\alpha}})\in\Xi$ arbitrary
$$Z(\alpha,{\bar{\alpha}}):=\frac{p^{\sum_{k=1}^4\alpha_k}(1-p)^{\sum_{k=1}^4{\bar{\alpha}}_k}}{(k_1\sum_{k=1}^4\alpha_k+k_2\sum_{k=1}^4{\bar{\alpha}}_k)^{1/2}}e^{-\frac{1}{2}\big[k_1\sum_{k=1}^4\alpha_k\theta_k^2+k_2\sum_{k=1}^4{\bar{\alpha}}_k\theta_k^2-\frac{(k_1\sum_{k=1}^4\alpha_k\theta_k+k_2\sum_{k=1}^4{\bar{\alpha}}_k\theta_k)^2}{k_1\sum_{k=1}^4\alpha_k+k_2\sum_{k=1}^4{\bar{\alpha}}_k}\big] },$$ 
which is the denominator in (\ref{40}). Next, by the change of variables
$$s=\frac{1}{\sqrt{k_1\sum_{k=1}^4\alpha_k+k_2\sum_{k=1}^4{\bar{\alpha}}_k+k_1-k_2}}\bigg[t+\frac{k_1\sum_{k=1}^4\alpha_k\theta_k+k_2\sum_{k=1}^4{\bar{\alpha}}_k\theta_k+(k_1-k_2)\theta_j}{\sqrt{k_1\sum_{k=1}^4\alpha_k+k_2\sum_{k=1}^4{\bar{\alpha}}_k+k_1-k_2}}\bigg],$$
in each of the sixteen ensuing Gaussian integrals of $\ey_{\nu_{x,\psi}}\left(-g''(\nabla_{j}\phi(x))\right)$, we obtain after integration
\begin{eqnarray*}
\lefteqn{\ey_{\nu_{x,\psi}}\left(-g''(\nabla_{j}\phi(x))\right)}\nonumber\\
&\le& p(1-p)\sqrt{2\pi k_1k_2}\,\,\,+\frac{\tilde{c}}{\sqrt{k_1-k_2}}\sum_{(\alpha,{\bar{\alpha}})\in\Xi} \frac{1}{Z}\frac{p(1-p)^{-1}(k_1-k_2)^2 p^{\sum_{k=1}^4\alpha_k}(1-p)^{\sum_{k=1}^4{\bar{\alpha}}_k}}{(k_1\sum_{k=1}^4\alpha_k+k_2\sum_{k=1}^4{\bar{\alpha}}_k+k_1-k_2)^{1/2}}\nonumber\\
&&\left(\frac{k_1\sum_{k=1}^4\alpha_k(\theta_k-\theta_j)+k_2\sum_{k=1}^4{\bar{\alpha}}_k(\theta_k-\theta_j)}{k_1\sum_{k=1}^4\alpha_k+k_2\sum_{k=1}^4{\bar{\alpha}}_k+k_1-k_2}\right)^2
e^{-(k_1-k_2)\theta_j^2-k_1\sum_{k=1}^4\alpha_k\theta_k^2-k_2\sum_{k=1}^4{\bar{\alpha}}_k\theta_k^2}\\
&& e^{\frac{\big(k_1\sum_{k=1}^4\alpha_k\theta_k+k_2\sum_{k=1}^4{\bar{\alpha}}_k\theta_k+(k_1-k_2)\theta_j\big)^2}{{k_1\sum_{k=1}^4\alpha_k+k_2\sum_{k=1}^4{\bar{\alpha}}_k+k_1-k_2}}}.
\end{eqnarray*}
Using inside each of the sixteen $(\alpha,{\bar{\alpha}})$ sums the lower bound $Z\ge Z(\alpha,{\bar{\alpha}})$, we get in the above
\begin{eqnarray}
\label{70}
\lefteqn{\ey_{\nu_{x,\psi}}\left(-g''(\nabla_{j}\phi(x))\right)}\nonumber\\
&\le& p(1-p)\sqrt{2\pi k_1k_2}\,\,\,+\tilde{c}p(1-p)^{-1}(k_1-k_2)^{3/2}\sum_{(\alpha,{\bar{\alpha}})\in\Xi} \frac{(k_1\sum_{k=1}^4\alpha_k+k_2\sum_{k=1}^4{\bar{\alpha}}_k)^{1/2}}{(k_1\sum_{k=1}^4\alpha_k+k_2\sum_{k=1}^4{\bar{\alpha}}_k+k_1-k_2)^{1/2}}\nonumber\\
&&\left(\frac{k_1\sum_{k=1}^4\alpha_k(\theta_k-\theta_j)+k_2\sum_{k=1}^4{\bar{\alpha}}_k(\theta_k-\theta_j)}{k_1\sum_{k=1}^4\alpha_k+k_2\sum_{k=1}^4{\bar{\alpha}}_k+k_1-k_2}\right)^2
e^{-(k_1-k_2)\theta_j^2+\frac{\big(k_1\sum_{k=1}^4\alpha_k\theta_k+k_2\sum_{k=1}^4{\bar{\alpha}}_k\theta_k+(k_1-k_2)\theta_j\big)^2}{{k_1\sum_{k=1}^4\alpha_k+k_2\sum_{k=1}^4{\bar{\alpha}}_k+k_1-k_2}}}\nonumber\\
&& e^{-\frac{(k_1\sum_{k=1}^4\alpha_k\theta_k+k_2\sum_{k=1}^4{\bar{\alpha}}_k\theta_k)^2}{k_1\sum_{k=1}^4\alpha_k+k_2\sum_{k=1}^4{\bar{\alpha}}_k}}.
\end{eqnarray}
Note now that 
\begin{eqnarray*}
\lefteqn{(k_1\sum_{k=1}^4\alpha_k\theta_k+k_2\sum_{k=1}^4{\bar{\alpha}}_k\theta_k+(k_1-k_2)\theta_j)^2}\\
&\le&(1+\lambda(\alpha,{\bar{\alpha}}))(k_1\sum_{k=1}^4\alpha_k\theta_k+k_2\sum_{k=1}^4{\bar{\alpha}}_k\theta_k)^2+\big(1+\frac{1}{\lambda(\alpha,{\bar{\alpha}})}\big)(k_1-k_2)^2\theta_j^2,
\end{eqnarray*}
where we choose $\lambda(\alpha,{\bar{\alpha}})>0$ such that 
$$\frac{1+\lambda(\alpha,{\bar{\alpha}})}{k_1\sum_{k=1}^4\alpha_k+k_2\sum_{k=1}^4{\bar{\alpha}}_k+k_1-k_2}<\frac{1}{k_1\sum_{k=1}^4\alpha_k+k_2\sum_{k=1}^4{\bar{\alpha}}_k}$$ and 
$$\frac{(k_1-k_2)(1+1/\lambda(\alpha,{\bar{\alpha}}))}{k_1\sum_{k=1}^4\alpha_k+k_2\sum_{k=1}^4{\bar{\alpha}}_k+k_1-k_2}<1.$$ 
Then
\begin{eqnarray*}
  \lefteqn{-(k_1-k_2)\theta_j^2+\frac{\big(k_1\sum_{k=1}^4\alpha_k\theta_k+k_2\sum_{k=1}^4{\bar{\alpha}}_k\theta_k+(k_1-k_2)\theta_j\big)^2}{{k_1\sum_{k=1}^4\alpha_k+k_2\sum_{k=1}^4{\bar{\alpha}}_k+k_1-k_2}}-\frac{(k_1\sum_{k=1}^4\alpha_k\theta_k+k_2\sum_{k=1}^4{\bar{\alpha}}_k\theta_k)^2}{k_1\sum_{k=1}^4\alpha_k+k_2\sum_{k=1}^4{\bar{\alpha}}_k}}\\
&\le& -(k_1-k_2)\theta_j^2+ \frac{\big(1+\frac{1}{\lambda(\alpha,{\bar{\alpha}})}\big)(k_1-k_2)^2\theta_j^2}{{k_1\sum_{k=1}^4\alpha_k+k_2\sum_{k=1}^4{\bar{\alpha}}_k+k_1-k_2}}\\
&&+ \frac{(1+\lambda(\alpha,{\bar{\alpha}}))(k_1\sum_{k=1}^4\alpha_k\theta_k+k_2\sum_{k=1}^4{\bar{\alpha}}_k\theta_k)^2}{{k_1\sum_{k=1}^4\alpha_k+k_2\sum_{k=1}^4{\bar{\alpha}}_k+k_1-k_2}}-\frac{(k_1\sum_{k=1}^4\alpha_k\theta_k+k_2\sum_{k=1}^4{\bar{\alpha}}_k\theta_k)^2}{k_1\sum_{k=1}^4\alpha_k+k_2\sum_{k=1}^4{\bar{\alpha}}_k}\\
&\le& -\epsilon_1(\alpha,{\bar{\alpha}},k_1,k_2)(k_1-k_2)\theta_j^2-\epsilon_2(\alpha,{\bar{\alpha}},k_1,k_2)\frac{(k_1\sum_{k=1}^4\alpha_k\theta_k+k_2\sum_{k=1}^4{\bar{\alpha}}_k\theta_k)^2}{k_1\sum_{k=1}^4\alpha_k+k_2\sum_{k=1}^4{\bar{\alpha}}_k},
\end{eqnarray*}
for some $\epsilon_1(\alpha,{\bar{\alpha}},k_1,k_2),\epsilon_2(\alpha,{\bar{\alpha}},k_1,k_2)>0$. Then (\ref{70}) becomes
\begin{eqnarray}
\label{estg''}
\lefteqn{\ey_{\nu_{x,\psi}}\left(-g''(\nabla_{j}\phi(x))\right)}\nonumber\\
&\le& p(1-p)\sqrt{2\pi k_1k_2}\,\,\,+\tilde{c}p(1-p)^{-1}(k_1-k_2)^{3/2}\sum_{(\alpha,{\bar{\alpha}})\in\Xi} \frac{(k_1\sum_{k=1}^4\alpha_k+k_2\sum_{k=1}^4{\bar{\alpha}}_k)^{1/2}}{(k_1\sum_{k=1}^4\alpha_k+k_2\sum_{k=1}^4{\bar{\alpha}}_k+k_1-k_2)^{1/2}}\nonumber\\
&&\left[2\left(\frac{k_1\sum_{k=1}^4\alpha_k\theta_k+k_2\sum_{k=1}^4{\bar{\alpha}}_k\theta_k}{k_1\sum_{k=1}^4\alpha_k+k_2\sum_{k=1}^4{\bar{\alpha}}_k+k_1-k_2}\right)^2+2\left(\frac{4(k_1+k_2)}{k_1-k_2}\right)^2\theta_j^2\right]
e^{-(k_1-k_2)\epsilon_1(\alpha,{\bar{\alpha}},k_1,k_2) \theta_j^2}\nonumber\\
&&e^{-\epsilon_2(\alpha,{\bar{\alpha}},k_1,k_2)\frac{(k_1\sum_{k=1}^4\alpha_k\theta_k+k_2\sum_{k=1}^4{\bar{\alpha}}_k\theta_k)^2}{k_1\sum_{k=1}^4\alpha_k+k_2\sum_{k=1}^4{\bar{\alpha}}_k}}
\le p(1-p)\sqrt{2\pi k_1k_2}+\epsilon_3 p(1-p)^{-1}\sqrt{k_1-k_2},
\end{eqnarray}
for some $\epsilon_3>0$ and where for the last inequality we have used $xe^{-x}<1$, with $x>0$, to bound the exponential part. Combining (\ref{lowboundbk}), (\ref{vx22}) and (\ref{estg''}), the conclusion follows.

\section*{Acknowledgment}
We thank Elliott Lieb for suggesting to us the use of the even/odd representation and Christof K\"ulske for pointing out to us that we can explicitly compute the 1-step iteration in example 3.2 (a). We also thank Nicolas Petrelis, Rongfeng Sun and an anonymous referee for many useful comments and suggestions, which greatly improved the presentation of the manuscript.

\end{document}